\documentclass{amsart}
\usepackage{enumerate,amssymb,amsmath,graphicx,ucs, tikz, mathtools}
\usepackage{wrapfig}
\usepackage{hyperref}
\usepackage{todonotes}
\usepackage[cp1251]{inputenc}
\usepackage{times,mathabx,amsfonts,amssymb,amsrefs,mathrsfs,tikz, xcolor}
\usepackage{todonotes}
\usetikzlibrary{decorations,decorations.pathmorphing}

\usetikzlibrary{matrix}
\usetikzlibrary{math, patterns}

\numberwithin{figure}{section}

\newcommand{\thistheoremname}{}

\newtheorem*{genericthm*}{\thistheoremname}
\newenvironment{namedthm*}[1]
{\renewcommand{\thistheoremname}{#1}%
	\begin{genericthm*}}
	{\end{genericthm*}}
\newtheorem{thm}{Theorem}[section]
\newtheorem{prop}[thm]{Proposition}
\newtheorem{lemma}[thm]{Lemma}

\newtheorem{claim}[thm]{Claim}
\newtheorem{corol}[thm]{Corollary}

\theoremstyle{definition}
\newtheorem{definition}[thm]{Definition}
\newtheorem{rem}[thm]{Remark}
\newtheorem*{thm*}{Theorem}
\newtheorem*{definition*}{Definition}
\newtheorem*{lemma*}{Lemma}
\newtheorem*{statement*}{Claim}
\newtheorem{example}[thm]{Example}
\newtheorem*{corol*}{Corollary}
\newtheorem*{question*}{Question}

\def\q#1.{{\bf #1.}}

\def\Z{\mathbb Z}

\newcommand{\MM}{\operatorname{M}}
\newcommand{\OR}{\operatorname{OR}}

\newcommand{\opt}{\operatorname{opt}}
\newcommand{\diam}{\operatorname{diam}}
\newcommand{\mins}{\operatorname{Br}}

\long\def\comment#1{}

\author{Anna Erschler}
\address{A.E.: C.N.R.S., \'{E}cole Normale Superieur, PSL Research University, France}
\email{anna.erschler@ens.fr}

\author{Ivan Mitrofanov}
\address{I.M.: C.N.R.S., \'{E}cole Normale Superieur, PSL Research University, France }
\email{phortim@yandex.ru }

\title{Assouad-Nagata 
dimension and gap for ordered metric spaces}

\date{\today}


\begin{document}
	\maketitle
	
	\begin{abstract}

		We prove that all spaces of finite Assouad-Nagata dimension admit a good order for Travelling Salesman Problem, and provide sufficient conditions under which  the converse is true. We formulate a conjectural characterisation of spaces of finite $AN$-dimension, which would  yield a gap statement for the efficiency of  orders on metric spaces. Under assumption of doubling, we prove a stronger gap phenomenon about all orders on a given metric space.
		
	\end{abstract}
	
	\section{Introduction}
	Given a metric space $(M,d)$, we consider  a finite subset $X$ of $M$.  We denote by  $l_{\opt}(X)$  the minimal length of a path which visits all points of $X$.
	Now assume that  $T$ 
	is a total order on $M$.
	For a finite subset $X \subset M$ we consider the restriction of the order
	$T$ on $X$, and enumerate 
	the points of $X$ accordingly:
	$$
	x_1 \leq_T x_2 \leq_T x_3 \leq_T \dots \leq_T x_{k}
	$$
	where $k = \#X$. Here and in the sequel  $\#X$ denotes the cardinality of the
	set $X$.
	We denote by $l_T(X)$ the length of the corresponding path
	$$
	l_T(X) := d(x_1, x_2) + d(x_2, x_3) + \dots + d(x_{k-1}, x_k).
	$$

	Given an ordered metric space $\MM=(M, d,T)$ containing at least two points
	and $k\ge 1$, we define the {\it order ratio function}
	$$
	\OR_{M,T}(k) := \sup_{X\subset M | 2 \le \# X \le k+1} \frac{l_T(X)}{l_{\opt}(X)}.
		$$
	See Definition \ref{def:OR} and Section \ref{sec:firstexamples} for more on this definition.

	We say that a metric space $M$ is {\it uniformly discrete}, if there exists $\delta>0$ such that for all pairs of points $x \ne y$ the
	distance between $x$ and $y$ is at least $\delta$.

	The travelling salesman problem aims to construct a cycle of minimal total length that visits each of $k$ given points.
	Bartholdi and Platzman introduced the idea to order all points of a metric space and then, given a $k$-point subset,  visit its points in the corresponding order
	\cite{bartholdiplatzman82}, \cite{bartholdiplatzman89}.
	Such approach is called {\it universal travelling salesman problem}.
	(One of motivations of  Bartholdi and Platzman was that this approach works
	fast for subsets  of a two-dimensional plane). 
	Their  argument implies a logarithmic upper bound for the function $\OR_{\mathbb{R}^2}(k)$. 
	
	For some metric spaces the function $\OR$ is even better. 
	Our result in \cite{ErschlerMitrofanov1} (Thm B) shows that the best possible situation, when $\OR(k)$ is bounded by above by a constant,  holds true for uniformly discrete 
	$\delta$-hyperbolic
	spaces.

	While in the original paper  of Bartholdi and Platzman it was suggested that such efficient behavior (with bounded ratio) holds for $\mathbb{Z}^d$,  it is known that it is not the case (unless $d=1$).
	The initial argument of Bartholdi and Platzman shows the existence of an order with  $\OR(k) \le {\rm Const} \ln k$. 
	Orders with at most logarithmic $\OR$ also exist on spaces with doubling property,
	by the result of Jia et al  \cite{jiadoubling}. 
	It seems  natural to conjecture that  this upper bound is optimal for $\mathbb{Z}^d$, $d\ge 2$. 
	A question of existence of an order on $\mathbb{Z}^2$ that violates  a logarithmic lower bound is asked in
	\cite{
Christodoulou}.
	Under additional assumption that the order is hierarchical  this logarithmic lower bound  (for a unit square, and hence for $\mathbb{Z}^2$) is proven in  \cites{eades2}, \cite{eades}, Thm.1, Section 3.3.

	Observe that among finitely generated groups only those
	that are virtually nilpotent  satisfy doubling property (see \ref{subsection:doubling} for definitions and background). We show that many groups of exponential growth also admit a logarithmic upper bound for an appropriate choice of the order, as we explain below.

	The notion of Assouad-Nagata dimension
	goes back to \cite{Nagata58},
	and the term {\it Nagata dimension} is used in  \cite{Assouad82}. 
	This notion provides a control of both local and global properties of the space.
	If the space is uniformly discrete, only large scales
	(and thus global properties) matter, and this notion coincides with {\it linearly  controlled metric dimension} (also called
	linearly controlled asymptotic dimension). For uniformly discrete
	spaces $AN$-dimension is a quasi-isometric invariant.
	We recall the definition of Assouad-Nagata dimension in Section \ref{section:finitedimension}. 
	Here we mention
	that the following spaces have
	finite $AN$-dimension.
	Spaces with doubling property \cite{LangSchlichenmaier},  wreath products of groups of linear growth with finite ones \cite{BrodskiyDydakLang}, polycyclic groups \cite{HigesPeng},
	trees (see e.g. \cite{roe}), 
	hyperbolic groups (and $\delta$-hyperbolic spaces that are "doubling in the small" \cite{LangSchlichenmaier}),
	and more generally, graphs and groups admitting quasi-isometric imbedding into finite products of trees (for imbedding of hyperbolic groups see \cite{BuyaloDranishnikovSchroeder}) such as Coxeter groups \cite{DranishnikovJanuszkiewicz}
	(and more generally for virtually special groups see \cite{HaglundWise08}).  
	Groups relatively hyperbolic with respect to groups of finite Assouad-Nagata dimension  also have finite Assouad-Nagata dimension \cite{Hume17}. 
	A not necessarily finitely presented $C'(1/6)$ small cancellation groups have $AN$-dimension at most $2$ \cite{Sledd1}.
	
	In fact, the argument of \cite{jiadoubling} for $\OR(k)$
	in the case of  doubling spaces can be adapted to prove a logarithmic upper bound more generally for spaces of finite Assouad-Nagata dimension. 
	In \cite{jiadoubling}
	the authors study the notion of $(\sigma, I)$-partitioning schemes, the existence of which
	can be shown to be equivalent to the finiteness of Assouad-Nagata dimension. 
	In the proof on the theorem below we argue directly in terms of $AN$-dimension.
	The second claim of the theorem below provides not only the asymptotic bound on $OR$, but also discusses the value  of the order breakpoint
		$\mins(M,T)$, which is  defined
		as the smallest integer $s$ such that  $\OR_{M,T}(s) < s$ (Definition \ref{def:mins}). For some basic examples of $\mins$, see also \cite{ErschlerMitrofanov1}, where this notion was introduced.
	
	\begin{namedthm*}{Theorem I}(=Thm \ref{thm:nagata}) 
		If $M$ is a
		metric space of finite Assouad-Nagata dimension $m$ with $m$-dimensional control function
		at most $Kr$, then there exists an order $T$ such that for all $k\ge 2$
		\begin{enumerate}
			\item   $\OR_{M,T}(k) \le C \ln k$, where a positive constant $C$ can be chosen depending on $m$ and $K$ only.
			
			\item $\mins(M,T) \leqslant 2m + 2$. Moreover, the elongations of snakes on $2m+3$ points are bounded by some constant depending on $m$ and $K$ only.
		\end{enumerate}
	\end{namedthm*}
	
	This theorem gives us upper estimation 
	on invariants of $M$:
	$\OR_M(k) = O(\log k)$ and
	$\mins(M) \leqslant 2m+2$,
	where
    $$
	\OR_{M}(k) := \inf_{T} \OR_{M,T}(k).
	$$

	"Snakes" mentioned in the second claim of the theorem are   order-increasing sequences of points which oscillate between neighborhoods of two points
	(We discuss the notion of snakes in more detail in Section \ref{sec:firstexamples}).
	
	In view of Lemma \ref{rem:goedel}, it is sufficient to prove the statement of Theorem \ref{thm:nagata}  for finite metric spaces. 
	We make a more general assumption that $M$ is uniformly discrete.
	We choose an appropriate constant $\lambda$, defined by the linearity constant $K$ for the control function
	(a possible choice is $\lambda= 4K$), and consider coverings from the definition
	of Assouad-Nagata dimension with $r= \lambda^n$, $n \in \mathbb{N}$. In Lemma \ref{lem:finitedimensionimpliesANfiltration} we modify this family of coverings to enforce a certain hierarchical structure
	on the sets of these coverings. This hierarchical structure  guarantees that the  sets of the coverings satisfy the assumption of Lemma \ref{le:convex}, which provides a sufficient
	condition for the existence of an order, such that 
	given  sets on some space are "convex"  with respect to this  order. We discuss this notion at the end of Section \ref{sec:firstexamples}.
	In Lemma \ref{lem:ANfiltrationCorollary}
	we show that such orders satisfy an upper bound for the order ratio function in the claim of the theorem.
	\bigskip
	
	The worst possible case for solving the universal travelling salesman problem are spaces with linear $\OR(k)$.
	An example of a sequence of finite graphs with linear $\OR(k)$ is constructed in Gorodezky et al  \cite{gorodezkyetal}, see also Bhalgat et al \cite{bhalgatetal}  who show that
	a sequence of Ramanujan graphs of large girth and of bounded diameter-by-girth ratio has 
	this property.
	The above mentioned paper considered both the question of the dependence of the number  of required points $k$ as well as the dependence on the cardinality of a finite graph $n$. As we have mentioned,  in this paper we consider the dependence on $k$, the question that makes both sense  for finite and infinite spaces.
	
	Since a result of Osajda
	\cite{osajda} allows to imbed subsequences of graphs with large girth into Cayley graphs,
	combining his result with that of \cite{gorodezkyetal}
	one can conclude that there exist groups with linear $\OR(k)$.
	While the above mentioned argument  uses both
	a remarkable construction of  Ramanujan graphs of large girth and a recent graphical small cancellation technique, we prove that 
	there is a large class of metric spaces and sequences of metric spaces 
	with infinite order breakpoint (and thus claiming that  $\OR(k)=k$, not only that this function is linear). Easy examples of groups of this kind can be obtained from Thm III. 
	
	Before stating this theorem, we formulate Thm II,  the first claim of which gives a sufficient spectral condition for infinite $\mins$. 
	Informally speaking, infinite $\mins$ means that whatever orders we choose on vertices of our graphs, there are arbitrary large  subsets of vertices on which this order is extremely far from optimal for  the travelling salesman problem.
	The second 
	claim of Thm II provides additional information about snakes for a sequence of expander graphs.
    In this theorem we use the notion of the order ratio function and the order breakpoint not only for a given space, but also for a sequence of graphs, which can be defined   analogously, see the beginning of Section~\ref{sec:weakexpansion}.

	\begin{namedthm*}{Theorem II}(=(1) and (3) of Thm \ref{thm:expanders})
		Let $\Gamma_i$ be a sequence of   finite graphs of degree $d_i\ge 3$
		on $n_i$ vertices.
		Let $T_i$ be an order on $\Gamma_i$.
		\begin{enumerate}
			\item  Assume that the normalized spectral gap
			$\delta_i= (\lambda^{\Gamma_i}_1 - \lambda^{\Gamma_i}_2 )  /d_i$ satisfies
			$$
			1/\delta_i  = o  \left( \frac{\log_{d_i}n_i}{\ln 
				\log_{d_i}n_i} \right).
			$$

			Then the order breakpoint of the sequence  $(\Gamma_i, T_i)$ is infinite. 
			
			\item For   a sequence of bounded degree expander graphs
			the following holds. If $d_i=d\ge 3$ and $\delta=\inf_i \delta_i >0$,
			then for each $k$ the graphs $(\Gamma_i, T_i)$ admit snakes on $k$ points of bounded width $C_k$ and of length at least $\log_{d-1}n_i -C'_k$, for some $C_k, C'_k>0$.
		\end{enumerate}
	\end{namedthm*}
	
	In Theorem II above we have formulated claims (1) and (3) of
	Thm \ref{thm:expanders}. In claim (2) of
	Thm \ref{thm:expanders}
	we will also give an estimate on possible length and width of a snake, in terms of
	a spectral gap of a graph.
	
	In general,  if the decay of $\delta_i$ is quick, one can not understand whether  $\mins$ is infinite or not, given a sequence of graphs
	and knowing  their cardinality $n_i$, degree  $d_i$ and spectral gap $\delta_i$,
	see Remark \ref{rem:nospectral}.
	For Claim (2) of the Theorem, the assumption of expansion can not be weakened if we want a criterion in terms of $n_i$, $d_i$ and $\delta_i$, see Remark \ref{rem:closetoexpandersnobs}. One can ask whether
	$o(\log_{d_i}^2(n_i))$
	(which is a sufficient condition for infinite $AN$-dimension of bounded degree graphs by \cite{humeetal}, see
	Remark \ref{rem:ANpoincare}) is
	sufficient for infiniteness of
	$\mins$.
	If 
	$1/\delta_i \sim (\log_{d_i}{n_i})^2$ then $\mins$ can be finite (see Remark \ref{rem:wreath}).
	Since a sequence of expander graphs can have a diameter close to $\log_{d-1}n_i$, the lower bound on the length of the snake in the claim (2)
	can not be significantly improved (see also Remark \ref{rem:sardari}).
	
	In Section \ref{sec:infinitegirth}
	we provide
	another sufficient condition of a different nature for infinite $\mins$.
	It can be deduced from Lusternik-Schnirelmann theorem that any order on an $\varepsilon$-net of a sphere $S^k$ admits snakes 
	on $k+2$ points, alternating between $\varepsilon$-neighborhoods of  antipodal points (see Lemma  \ref{lem:sphere}). Combining it with the control 
	of order ratio function for weak imbeddings of cubes we get
	
	\begin{namedthm*}{Theorem III}(= Corollary \ref{cor:cubes})
		If a metric space $M$ weakly contains  arbitrarily large cubes of dimension $d$,
		then for any order $T$ on $M$ it holds
		$$
		\OR_{M,T}(d) = d.
		$$
		
		\noindent In particular, if a metric space $M$ weakly contains a sequence of arbitrarily large cubes,
		then for any order $T$ on $M$ the order breakpoint of $(M, T)$ is infinite.
	\end{namedthm*}
	
	We recall again that the informal meaning of infinite order breakpoint is that whatever order we choose on $M$, this order behaves extremely bad (for the travelling salesman problem) on some arbitrary large subsets.
	
	In Section \ref{sec:infinitegirth} we will define for metric spaces the property of containing a sequence of arbitrarily large cubes.  
	Here  we mention
	that the class of such spaces includes spaces admitting uniform imbeddings of $\mathbb{Z}^d$ (or $\mathbb{Z}_+^d$)
	for all $d$.
	In particular, this condition holds for any  finitely generated group $G$ that
	contains the direct sum $\mathbb{Z}^\infty$ as a subgroup. 
	These include many classes of solvable groups and many other classes of amenable groups.
	Here we mention that Grigorchuk groups (and moreover all known constructions of groups of intermediate growth) admit uniform
	imbeddings of $\mathbb{Z}_+^d$, for all $d$.
	These also include many examples of  not amenable groups, as well as
	some groups such as Thompson group (where famous amenability question remains open).
	Further examples of spaces that weakly contain sequences of arbitrarily large cubes  are  
	$\mathbb{Z}^2\wr \mathbb{Z}/2\mathbb{Z}$ and more generally $B \wr A$, where $B$ is an infinite group of not linear growth and $A$ is any finite or infinite group containing at least two elements (this statement is inspired by the argument in \cite{BrodskiyDydakLang},
	see Lemma \ref{lem:wcubeswr} where we  study weak imbeddings of cubes in  wreath products). We also mention that  imbeddings of cubes appear naturally as lower estimates for $AN$-dimension in various classes of groups and spaces, see \cites{higes, Sledd}.

	In view of  Theorems II and III mentioned above we ask
	
	\begin{question*}
		Let $M$ be a metric space
		of infinite Assouad-Nagata dimension. Is it true that the order 
		breakpoint of $M$ is infinite? 
	\end{question*}
	We recall that there are various known examples of groups of finite asymptotic dimension which have infinite $AN$-dimension (see  Nowak \cite{nowak}, Brodskyi Dydak Lang \cite{BrodskiyDydakLang}).
	As we have already mentioned, some of them satisfy the assumption of our Theorem III.
	In view of Thm I, if the answer to the question is positive, this would provide an equivalent
	characterization of spaces  of finite $AN$-dimension. 
	Taking in account the above mentioned examples, for a question of a possible characterization it is essential to speak about $AN$-dimension and not
	about asymptotic dimension. 
	
	Observe also that if the answer to the above mentioned question is positive, in view of Theorem \ref{thm:nagata}
	this would provide a positive answer to the following
	
	\begin{question*}[Gap problem for existence of orders]
		Let $M$ be a metric space. Is it true that either for any order $T$ on $M$ and all $k\ge 1$ it holds
		$\OR_{M,T}(k) = k$
		or  there exists an order $T$ such that for all $k\ge 2$ it holds $
		\OR_{M,T}(k) \le \rm{Const} \ln k$?
	\end{question*}
	
  	Among classes of spaces where we obtain a positive answer to this question are Cayley graphs of wreath products of groups (see Corollary \ref{cor:wreathgap}).
	
	Given a metric space, one can formulate a stronger Gap problem, which describes behavior of all orders (rather
	then searches an order on the space). Our next result below solves this problem for
	spaces with doubling property.

	As we have already mentioned, the argument of Bartholdi and Platzman for Euclidean plane (and generalizations of their argument for the spaces with doubling property)
	provides orders with logarithmic
	upper bound for the function $\OR(k)$.
	This is in contrast with the lexicographic order on $\mathbb{R}^2$ ($(x_1, y_1) < (x_2, y_2)$ if $x_1 <x_2$, or $x_1=x_2$ and $y_1<y_2$), where it is easy to see that $\OR(s)=s$ for all $s$. 
	Our theorem below shows that any order on a space with doubling property (in particular, any order on $\mathbb{R}^d$) satisfies the same dichotomy:

	\begin{namedthm*}{Theorem IV}[Gap for order ratio functions on spaces with doubling property] (=Thm \ref{thm:gap})
		Let $M$ be a metric space with doubling property and $T$ be an order on $M$. Then
		either for all $s$ it holds
		$$
		\OR_{M,T}(s)= s
		$$
		or there exists $C$ (depending only on the doubling constant  of $M$,  $s$ and $\varepsilon$ such that $\OR_{M,T}(s)\le s- \varepsilon$) such that for all $k \ge 1$
		$$
		\OR_{M,T}(k) \le  C \ln k
		$$
	\end{namedthm*}
	
	{\bf Acknowledgements. } 
	We would like to thank  David Hume for discussions about $AN$ dimension,  Tatiana Nagnibeda for explanations and references on spectra of finite graphs;
	Karim Adiprasito for helpful conversations. We are grateful to the referees for useful comments.
	This project has received funding from the European Research Council (ERC) under
	the European Union's Horizon 2020 
	research and innovation program (grant agreement
	No.725773).
	The work of the second named author is also supported by Russian Science Foundation (grant no.17-11-01377).

	\section{Preliminaries and basic properties}\label{sec:firstexamples}

		Below we recall some basic properties of the order ratio function and order breakpoint, discussed  in \cite{ErschlerMitrofanov1}.

	\begin{definition}\label{def:OR} [Order ratio function]
		Given an ordered metric space $\MM=(M, d,T)$ containing at least two points
		and $k\ge 1$, we define the {\it order ratio function}
		$$
		\OR_{M,T}(k) := \sup_{X\subset M | 2 \le \# X \le k+1} \frac{l_T(X)}{l_{\opt}(X)}.
		$$

	If $M$ consists of a single point, then the supremum in the definition above is taken over an empty set. 
	We use in this case the convention that $\OR_{M,T}(k)=1$ for all $k\ge 1$.

	Given an (unordered) metric space $(M,d)$, we also  define the {\it order ratio function} as 
	$$
	\OR_{M}(k) := \inf_{T} \OR_{M,T}(k).
	$$
	\end{definition}
	
	Denote by $L$ the diameter of $X$. 
	It is clear that $L \leqslant l_{\opt}(X) \leqslant l_T(X)$
	and $l_T(X) \leqslant L(\#X - 1)$. 
	Hence $1 \leq \OR_{M,T}(k) \leq k$ for any $M$, $T$ and $k \ge 1$.

	\begin{definition} \label{def:mins} [Order breakpoint]
		Let $M$ be a metric space, containing at least two points,
		and let $T$ be an order on $M$.  
		We say that the order  
		breakpoint 
		$\mins(M,T) = s$ if $s$ is the smallest integer such that $\OR_{M,T}(s) < s$. 
		If such $s$ does not exist, we say that $\mins(M,T)=\infty$.
		A possible convention for a one-point space $M$ is to define $\mins(M,T)=1$.
	
	Given an (unordered)
	metric space $M$, we define $\mins(M)$
	as the minimum over all orders $T$ on $M$:
	$$
	\mins(M) = \min_{T}
	\mins(M,T).
	$$
	\end{definition}	
It is clear that for any $k>\mins(M,T)$ it holds $\OR_{M,T}(k) < k$.	
	
\noindent The property $\mins_{M,T} \geq k$ means that on some $k$-point subsets of $M$ the order $T$ behaves extremely non-optimal as an universal order for the universal travelling salesman problem.

 We will give an idea of such non-optimal subsets.
Consider a sequence of  $s$  points: $x_1 <_T x_2 <_T \dots <_T x_s$, let $a$ be the diameter of this set and let $b$ be the maximal  distance $d(x_i, x_j)$ where $i$ and $j$ are of the same parity. 
Then $a$ is called the {\it length} of the {\it snake} $(x_i)$, $b$ is called its {\it width} and the ratio $a/b$ is the {\it elongation} of the snake $(x_i)$. We say that elongation is $\infty$ if $b = 0$.

If a set $X$ consists of $k+1$ points and
the value of $l_T(X)$ is close to $kl_{\opt}(X)$, then it can be shown that any two points with indices of the same parity are relatively close to each other, compared with the diameter of $X$. 
Hence, one can show

\begin{lemma}[\cite{ErschlerMitrofanov1}, Lemma 2.11]
Let $(M,T)$ be an ordered metric space. Given $s \geqslant 2$, we have $\OR_{M,T}(s) = s$ if and only if in $(M,T)$ there exist snakes  of arbitrary large elongation on $s+1$ points.
\end{lemma}

If we allow a metric to assume the value $\infty$, then a natural way to define $\OR$ is to consider the ratio $l_T(X)/l_{\opt}(X)$ for  subsets $X$ with finite diameter. In this setting we can speak about $\OR$ and $\mins$ of disjoint unions of metric spaces.

\begin{wrapfigure}{l}{0.2\textwidth}
\includegraphics[width=0.2\textwidth]{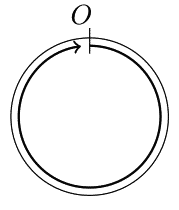}
\end{wrapfigure}
We mention a simple example (which we discussed already in \cite{ErschlerMitrofanov1}, Lemma 3.1) when the metric space $M$
is a circle $S^1$ with its inner metric.
It is not difficult to show, that $\OR_{S^1}(k)=2$ for all $k \geqslant 2$.
A natural order (a clockwise order)
provides the estimation $\OR_{S^1}(k)\leq 2$, and the lower estimation follows from

\begin{lemma}\label{le:examplecircle}
    For any order $T$ of $S^1$ and any $\varepsilon > 0$ there exist two antipodal points $x,y$ of the circle and a snake $z_1<_T z_2 <_T z_3$ on $3$ points such that $z_1$ and $z_3$ belong to  the $\varepsilon$-neighborhood of $x$ and $z_2$ is in the $\varepsilon$-neighborhood of $y$.
\end{lemma}
\begin{proof}
We call a point $x\in S^1$ {\it $T$-small} if $x<_T \Bar{x}$, where $\Bar{x}$ denotes the antipodal point to $x$.
Otherwise we call it {\it $T$-large}.
Since $(S^1,T)$ contains both $T$-small and $T$-large points, we can find two points $s_1,s_2$ of different types that are $\varepsilon$-close to each other.
Assume $s_1 <_T \Bar{s_1}$ and 
$s_2 >_T \Bar{s_2}$.
If $\Bar{s_1}<_T s_2$, we take the snake $s_1 <_T \Bar{s_1} <_T s_2$. 
Otherwise, we take the snake 
$\Bar{s_2} <_T s_2 <_T \Bar{s_1}$.
\end{proof}

We recall the notion of quasi-isometric imbedding.
Given metric spaces $N$ and $M$, a map $\alpha$  from $N$ to $M$ is a {\it  quasi-isometric imbedding} if
		there exist $C_1$, $C_2$ such that for any $x_1, x_2 \in N$ it holds
		\[
		\frac{1}{C_1}(d_N(x_1, x_2)) -C_2  \le   d_M(\alpha(x_1), \alpha(x_2)) \le C_1(d_N(x_1, x_2)) +C_2
		\]
		
		If $M$ is at bounded distance from $\alpha(N)$, this map is called a {\it quasi-isometry}, and the spaces $X$ and $Y$ are called {\it quasi-isometric.}
		If $\alpha$ is bijective and $C_2=0$, then the spaces are said to be  {\it bi-Lipschitz equivalent}.
	
We also recall	a  weaker condition of
uniform imbeddings.
	  Given  metric spaces $N$ and $M$, a map $\alpha$  from $N$ to $M$ is an
		{\it uniform imbedding} (also called {\it coarse imbedding})
		if there exist two non-decreasing functions $\rho_1$, $\rho_2:$ $[0, +\infty) \to [0, + \infty)$,
		with $\lim_{r\to \infty} \rho_1(r) = \infty$,  such that
		\[
		\rho_1 (d_X(x_1, x_2)) \leqslant d_M(\alpha(x_1), \alpha(x_2)) \le \rho_2 (d_X(x_1, x_2)).
		\]

		Given $\varepsilon, \delta > 0$, a subset $U$ in a metric set $M$ is
		said to be an $(\varepsilon, \delta)${\it -net} if any point of $M$ is at distance at most
		$\varepsilon$ from $U$ and the distance between  any two distinct points of $U$ is at least $\delta$.
It is easy to show that for any metric space $M$ and any $\varepsilon > 0$ there exists an 	$(\varepsilon, \varepsilon)$-net of $M$.

\begin{definition}
    We say that two functions $f_1,f_2(r) :\mathbb{N} \to \mathbb{R}$ are equivalent up to a multiplicative constant if  there exists $K_1,K_2>0$ such that $f_1(r) \le K_1 f_2(r)$ and $f_2(r) \le K_1 f_1(r)$  for all $r\ge 1$.
\end{definition}

It is easy to observe (\cite{ErschlerMitrofanov1}, Lemma 2.9) that if two metric spaces $M$ and $N$ are quasi-isometric, $M'$ is an $(\varepsilon, \delta)$-net of $M$ and $N'$ is an $(\varepsilon, \delta)$-net of $N$, then the order ratio functions $\OR_{M'}$ and $\OR_{N'}$ are equivalent up to a multiplicative constant.
In particular, the asymptotic class of $\OR$ is a quasi-isometric invariant of uniformly discrete metric spaces. Moreover, if a metric space $N$ is quasi-isometrically imbedded into a space $M$, then $\OR_M$ is asymptotically larger or equal to $\OR_N$.
We mention in this context
that given (a not necessary injective) map $\phi: N \to M$ and an order $T_M$ on $M$, one can construct an order $T_N$ on $N$ such that $\phi(x)<_{T_M}\phi(y)$ always implies $x<_{T_N}y$. 
We call any such order a {\it pullback of} $T_M$.

While the asymptotic class of the order ratio function is invariant under quasi-isometries, specific values of $\OR$ can change under quasi-isometries, as can be easily seen already from finite space examples.
But for any $s$ the equality $\OR(s) = s$ is preserved under quasi-isometries of uniformly discrete spaces (\cite{ErschlerMitrofanov1} Lemma 2.12).
To show this one can consider  a pullback and argue that snakes of large elongation map under quasi-isometries to snakes of large elongation.
Moreover, an uniformly discrete metric space can not be quasi-isometrically imbedded to a metric space with smaller $\mins$.

From the definition it is clear that if $M$ has at least two points then $\mins(M,T) \ge 2$ for any order $T$. 
A result of M.~Kapovich \cite{kapovich} about $(R,\varepsilon)$-tripods
can be used to characterise uniformly discrete metric spaces
with $\mins \le 2$ (see
Lemma 4.1 in  \cite{ErschlerMitrofanov1}): such spaces are either bounded or quasi-isometric to a ray or to a line.
It is not difficult to see that virtually free groups have $\mins \le 3$, and
Theorem $A$ of the above mentioned  paper shows the converse.
Let $G$ be a finitely generated group. 
Then $G$ admits an order $T$ with  $\mins(G,T) \le 3$	if and only if $G$ is virtually free.

While the property $\OR_M(s) = s$ is not preserved under uniform mappings of uniformly discrete spaces, 
there is a similar  property which is inherited by such imbeddings. Below  we recall this property, studied in \cite{ErschlerMitrofanov1}.
We say that
an ordered metric space $(M, T)$ admits {\it a sequence of  snakes of bounded width} on $s$ points
if there is  a sequence of snakes  on $s$ points of uniformly bounded width and with diameters tending to infinity.
It is not difficult to see that
if $N$ admits an uniform imbedding   to a metric space $M$ and 
if $N$ admits a sequence of  snakes of bounded width on $s$ points for any order on this space, 
then the same happens for $M$;
and in particular $\OR_M(s-1) = s-1$ (\cite{ErschlerMitrofanov1}, Lemma 2.14).

Using a compactness argument one obtains the following.
	
 	\begin{lemma}\label{rem:goedel}[\cite{ErschlerMitrofanov1}, Lemma 2.15]
		Let $M$ be a metric space. Consider a  function $F:\mathbb{Z}_+ \to \mathbb{R}_+$ and assume that for any finite subset $M' \subset M$ there exists an order $T'$ satisfying
		$\OR_{M', T'} (k) \le F(k)$ for all $k\ge 1$.
		Then there exists an order $T$ on $M$ satisfying 
		$\OR_{M, T} (k) \le F(k)$ for all $k\ge 1$.
	\end{lemma}
Using an induction on the cardinilty of  $\mathcal{A}$ in case when this family is finite  and a compactness argument  for the general case, one also obtains
	\begin{lemma} \label{le:convex}[\cite{ErschlerMitrofanov1}, Lemma 2.17] 
		Let  $\mathcal{A}$ be a family of subsets of  $M$ such that for
		any two sets $A_1, A_2 \in \mathcal{A}$ either their intersection is empty, 
		or one of the sets is contained in the other one.
		Then there exists an order $T$ on $M$ such that for any set $V\in \mathcal{A}$ and any three points $x,y,z\in M$ such that $x <_T y <_T z$, the condition $x, z\in V$ implies $y\in V$.
	\end{lemma}
	
	We say that the subset $V$ of $M$ is {\it convex with respect to the order} $T$ if $V$ and $T$ satisfy the property described in Lemma \ref{le:convex}.
	
This elementary lemma can be used
to construct orders with given properties in various spaces. We will use it several times in the proof of Theorem I about spaces of finite AN dimension and of Theorem IV about spaces with doubling property.

\section{Gap for the order ratio functions. Nilpotent groups and spaces with doubling property.} \label{sec:doubling}
	
	We have mentioned that a crucial observation for using orders for the travelling salesman  problem
	is the result of Bartholdi and Platzman, who have shown existence of orders (related to space-filling curves) such that for any finite set $X$ in the plane the order provides a tour with length at most $\rm{Const} (\ln \#X)$ times longer than the optimal one. 
	
	Consider two orders on the square $K = [0;1]^2$.
	The first order $T_{BP}$ is  the one considered by  Bartholdi and Platzman in \cite{bartholdiplatzman82} and their estimate
	(combined with obtained later lower bound in
	\cite{bertsimasgrigni})
	implies  that $\OR_{K,T_{BP}}(k)\sim \log k$.
	
	The second order 
	$T_{lex}$ is the lexicographical order. 
	We define this order by saying 
	$(x_1,y_1) <_{T_{lex}} (x_2, y_2)$ if
	$x_1 < x_2$ or $x_1 = x_2$ and $y_1 < y_2$. It is  clear that $(K,T_{lex})$ 
	contains snakes of diameter $1$ and arbitrary small width, and hence
	$\OR_{K,T_{lex}}(k) = k$ for all $k$.
	
	\medskip
	
	\newcommand{\bpone}[3]
	{\begin{scope}[shift={(#1,#2)},rotate=#3]
			\draw[ thick, rounded corners, ->] (0.3,0.1)--(0.9,0.1)--
			(0.9,0.9) -- (1.1,0.9) -- (1.1,0.1) -- (1.7,0.1);
	\end{scope}}
	
	\newcommand{\bptwo}[3]
	{\begin{scope}[shift={(#1,#2)},rotate=#3]
			\bpone{0}{0}{0}
			\bpone{2}{0}{90}
			\bpone{2}{2}{270}
			\bpone{2}{0}{0}
			\draw[ thick, rounded corners](1.7,0.1) -- 
			(1.9, 0.1) -- (1.9, 0.3);
			\draw[ thick, rounded corners](2.3,0.1) -- 
			(2.1, 0.1) -- (2.1, 0.3);
			\draw[ thick, rounded corners](1.9,1.7) -- 
			(1.9, 1.9) -- (2.1, 1.9) -- (2.1, 1.7);
			
	\end{scope}}
	
	\newcommand{\bpthree}[3]
	{\begin{scope}[shift={(#1,#2)},rotate=#3]
			\bptwo{0}{0}{0}
			\bptwo{4}{0}{90}
			\bptwo{4}{4}{270}
			\bptwo{4}{0}{0}
			\draw[ thick, rounded corners](3.7,0.1) -- 
			(3.9, 0.1) -- (3.9, 0.3);
			\draw[ thick, rounded corners](4.3,0.1) -- 
			(4.1, 0.1) -- (4.1, 0.3);
			\draw[ thick, rounded corners](3.9,3.7) -- 
			(3.9, 3.9) -- (4.1, 3.9) -- (4.1, 3.7);
			
	\end{scope}}
	
	\resizebox{13cm}{6cm}{
		\begin{tikzpicture}[scale = 1.5]
		\bpthree{0}{0}{0}
		\bpthree{8}{0}{90}
		\bpthree{8}{8}{180}
		\bpthree{0}{8}{-90}
		\draw[ thick, rounded corners](7.7,0.1) -- 
		(7.9, 0.1) -- (7.9, 0.3);
		\draw[ thick, rounded corners](7.9,7.7) -- 
		(7.9, 7.9) -- (7.7, 7.9);
		\draw[ thick, rounded corners](0.3,7.9) -- 
		(0.1, 7.9) -- (0.1, 7.3);
		
		\begin{scope}[shift = {(10,0)}]
		\foreach \x in {0,...,16}{
			\draw[thick, ->] (\x/2,4) -- (\x/2,8) -- (\x/2+0.5,0) -- (\x/2 + 0.5, 4);	
		}
		\end{scope}
		\end{tikzpicture}
	}
	
	\medskip

	Our result below will prove that neither on $K$ nor on any other space with doubling property it is possible to
	construct an order $T$ such that $\OR_{K,T}(k)$ is sub-linear but grows faster than a logarithmic function in $k$.
	
	\subsection{Spaces with doubling property}\label{subsection:doubling}
	
	\begin{definition}
		A metric space $(M, d)$ is said to satisfy {\it doubling property} if 
		there exists a doubling constant $D>0$ such that 
		for any $x\in M$ and $r>0$ the ball $B(x,2r)$ of radius $2r$ can be covered by at most $D$ balls of radius $r$. 
	\end{definition}
	Such spaces are also called doubling metric spaces.
	
	\noindent In subsection \ref{subs:gap_proof} we will prove the following
	
	\begin{thm}[Gap for order ratio functions on doubling spaces] \label{thm:gap}
		Let $X$ be a doubling metric space and $T$ be an order on $X$. Then
		either for all $k$ it holds
		$$
		\OR_{X,T}(k)= k
		$$
		or there exists $C$, depending only on the doubling constant of $X$,
		on the order breakpoint $s$ of $(X,T)$ and on $\varepsilon$ such that $\OR_{X,T}(s)\le s- \varepsilon$, such that for all $k\ge 2$
		$$
		\OR_{X,T}(k) \le  C \ln k.
		$$
	\end{thm}
	
	\noindent Spaces $\mathbb{R}^d$ and their subsets are examples of doubling spaces. It is easy to see that
	any group of polynomial growth has doubling property for infinitely many $r$.
	It is also known in fact that any group of polynomial growth has doubling property for all $r$, but known proofs use that these groups are virtually nilpotent (by
	Polynomial Growth theorem \cite{gromovpolynomial}).
	
	Below we recall the definition and state these remarks more precisely.
	\begin{definition} A {\it growth function} of a finitely generated group $G$ with respect to a finite generating set $S$ is the number of elements of $G$ of the word length at most $n$
		$$
		v_{G,S}(n) = \# B_{G,S}(n) = \# \{ g:  l_{G,S} (g) \le n\}
		$$
		
	\end{definition}
	
	\noindent Suppose that a group $G$ has polynomial growth, that is $v_{G,S} (n) \le C_S n^d$ for some (and hence for all)
	generating set $S$. It is clear directly from the definition  that there exists infinitely many $n$ such that $v_{G,S} (2 n ) \le C v_{G,S}(n)$ for some constant $C$.
	Using Polynomial Growth Theorem, we can  claim that $G$ is virtually nilpotent, and in this case
	in holds $ C_1 n^d  \le v_{G,S}(n) \le C_2 n^d $ for some integer $d$ (and some positive constants $C_1$ and $C_2$, depending on $S$; for this basic fact about growth see e.g. Chapter 4 of
	\cite{MannHowgroupsgrow}
	). In particular, there exists a positive constant $C$ depending on $S$  such that
	$$
	v_{G,S} (2 n ) \le C v_{G,S}(n)
	$$

	\noindent
	for all $n$.

	\begin{rem} Let $G$ be a finitely generated group, $R>0$. 
		Let $U$ be a $(R,R)$-net of the ball of radius $2R$  in $G$. Consider the balls of radius $R$ centered at points of the nets. 
		The family of these balls covers the ball of radius $2R$, and the number of the balls is at most $v(2.5 R)/v(R/2)$.
		In particular, if $v(2.5R)/v(R/2) \le v(2^3 R/2)/ v(R/2) \le C^3$, then $G$ has doubling property on the scale $R$ with constant $C^3$. 
	\end{rem}

	\subsection{Proof of Theorem \ref{thm:gap}}\label{subs:gap_proof}
	
	We  assume that there exists $s_0$ such that $\OR_{M,T}(s_0
	) = s_0 - \varepsilon$ for some $\varepsilon > 0$.
	Observe that
	there exists a number $\lambda$, depending only on $s_0$ and $\varepsilon$,
	such that any snake on $s_0+1$ points in $(X,T)$ has elongation (the ratio of its diameter and its  width) at most $\lambda$.
	
	\begin{lemma}\label{le:gap1}
		Let $M$ be a doubling metric space with a doubling constant $D$, let $T$ be an order on $M$ and let $\lambda>1$ be such that all snakes on $s_0+1$ points in $(M,T)$ have elongations at most $\lambda$.
		Then there exists a constant $N(D, s_0, \lambda)$
		such that 
		whatever  $x\in M$ and $R > 0$  we choose, there are no sequences 
		$x_1 <_T x_2 <_T \dots <_T x_{2N}$ satisfying
		\begin{enumerate}
			\item $x_i \in B(x,4R)$ for all $i$
			\item $d(x_{2i-1},x_{2i}) \geqslant R$ for all $i$.
		\end{enumerate}
		
	\end{lemma}
	
	\begin{proof}
		Arguing by induction on $k$ and using the definition of doubling spaces, we observe that for any $k$, any $x$ and
		any $R>0$ the ball $B(x,4R)$ can be covered by $D^{k+3}$ balls of
		radius at most $2^{-k-1}R$, and thus of 
		diameter not greater than $2^{-k}R$. 
		Choose $k$ such that $2^k > \lambda$.
		Put $N= s_0(D^{k+3})^2 + 1$. 
		Assume that there exist points
		$x_1 <_T x_2 <_T \dots <_T x_{2N}$ satisfying properties 
		(1) and (2) in the formulation
		of the Lemma.
		Cover  the ball
		$B(x,4R)$ by $m \leqslant D^{k+3}$ balls
		$B_1, \dots, B_m$ of diameters at most  $4R/2^{k+2} < \frac{R}{\lambda}$.
		We have $N$ pairs $(x_i,x_{i+1})$ of points at distance at least $R$. 
		Each pair $(x_i, x_{i+1})$ has $x_i$ in $B_j$ and $x_{i+1}$ in $B_l$ for some $j \neq l$.
		
		For given $j$ and $l$, 
		it is clear that there is at most  $s_0$ such pairs are in  $B_j \cup B_l$, since
		otherwise we find a snake on $s_0+1$ points with diameter at least $R$ and width smaller than $R/\lambda$. This is in a
		contradiction with the assumption of the lemma.
	\end{proof}
	
	\begin{claim}\label{le:dividepath}
		Let $M$ be a metric space and let $(x_0, x_1, x_2,\dots,x_n)$ be a sequence of points such that 
		\[\sum_{i=0}^{n-1} d(x_i, x_{i+1}) = L.\]
		Let $L = a_1 + a_2 + \dots + a_m$, all $a_i > 0$.
		Then the set $\{x_i\}$ can be partitioned into $m$ subsets 
		$A_1,\dots, A_m$ such that $\diam(A_j) \leqslant a_j$ for all $j$.
	\end{claim}
	
	\begin{proof}
		Induction on $m$. The base of induction: if $m=1$, then  the claim of the lemma is straightforward. 
		Induction step: find the maximal $l\ge 0$ such that 
		$\sum_{i=0}^{l-1} d(x_i, x_{i+1}) \leqslant a_1$ and put $A_1:= \{x_0, x_1,\dots,x_l\}$.
		Then $\sum_{i=l}^{n-1} d(x_i, x_{i+1}) \leqslant a_2 + \dots + a_m$ and $\diam(A_1) \leqslant a_1$.
	\end{proof}
	
	\begin{lemma}\label{le:gap3}
		Let $M$ be a doubling metric space with a doubling constant $D$ and let $T$ be an order on $M$ such that all snakes on $s_0+1$ points in $(M,T)$ have elongation at most $\lambda$.
		Then there exists a constant $N(D, s_0, \lambda)$ such that the following holds.
		
		Take a sequence of points $(x_i)$ in $M$, $i: 1 \le i \le m$  satisfying
		$$
		x_1 <_T x_2 <_T \dots <_T x_m.
		$$
		Assume that the minimal length of a path visiting  these points is $L$. Then  for all $k\ge 0$ it holds  
		\[
		\# \left\{ i \in (1,\dots, m-1) \Big| \frac{L}{2^{k+1}} \le d(x_i, x_{i+1}) \le \frac{L}{2^{k}}\right\} \leqslant N 2^k.
		\]
	\end{lemma}
	
	\begin{proof}
		Applying Claim \ref{le:dividepath}
		for $m=2^k$, $a_1=a_2 = \dots =a_{m}=L/2^k$ and the reordering of $\{x_i\}$ to minimal length 
		we conclude that there exist $2^{k}$ points $y_1,\dots y_{2^{k}}$ such that
		\[\{x_i\}_{i=1}^m \subset B\left(y_1, \frac{L}{2^{k}}\right) \bigcup \dots \bigcup B\left(y_{2^{k}}, \frac{L}{2^{k}}\right).
		\]
		Then for each $i$ such that $d(x_i,x_{i+1}) \leqslant \frac{L}{2^{k}}$ we can find $y_j$ such that 
		\[
		x_i, x_{i+1} \in B\left(y_j, \frac{L}{2^{k-1}}\right).
		\]
		Put $R:=2^{-k-1}L$ and apply Lemma \ref{le:gap1} for each $y_i$ and $R$.
		We can take $N(D, s_0, \lambda)$ twice as big as in Lemma \ref{le:gap1}.
	\end{proof}
	
	Now we are ready to complete the proof of Theorem \ref{thm:gap}.
	Consider a finite subset $X\subset M$ of cardinality $n$ and enumerate its points $x_1 <_T \dots <_T x_n$. Let $2^{k-1} < n \leqslant  2^k$.
	As we have already mentioned, there exists $\lambda$, depending on  $s_0$ and $\varepsilon>0$ satisfying $\OR(s_0)= s_0 -\varepsilon$, such that all snakes on $(s_0 + 1)$ points in $(M,T)$ have elongation at most $\lambda$.
	Consider $N=N(D, s_0,\lambda)$ from the claim of  Lemma \ref{le:gap3}.
	Let us show that
	\[
	l_{T}(X) = d(x_1,x_2) + \dots + d(x_{n-1}, x_n) < N l_{\opt}(Y)k.
	\]
	Denote $l_{\opt}(X)$ by $L$.
	
	We consider the following $n-1$ numbers: the distances $d(x_i,x_{i+1})$ for $1\leqslant i < n$. For any $i, i = 0,\dots, k-1$
	let $D_i$ be the number of distances between $L/2^{i+1}$ and $L/2^i$, and denote by $D_k$ the number of distances that are smaller than $2^{-k}L$. 
	Observe that each distance $d(x_i,x_{i+1})$ is at  most $L$.

	Note that 
	$D_i \leq N2^i$ for any $i$. 
	For $D_k$ it is clear  because $n < N2^k$, and for all other $D_i$ it follows from Lemma \ref{le:gap3}.
	Hence, the sum of numbers in each group is not greater than $NL$, and the total sum $l_T(X)$ is not greater than $NLk = Nl_{\opt}(X) \lceil \log_2(n) \rceil$.
	
	\subsection{Examples of orders with finite order breakpoint.}
	\label{subsec:euclid}
	
	Consider some  order $T$
	on $X = \mathbb{R}^d$ and suppose 
	that we want  to prove that $\OR_{X,T}(k) = O(\log(k))$.
	From Theorem \ref{thm:gap} it follows that its enough to find $s$ such that $\OR(s) < s$, i.e. such that
	the elongations of snakes on $s+1$ points in $(X,T)$ are bounded.
	
	Here we give an example of a simple order with this property, and together with Theorem \ref{thm:gap} this will give us an alternative proof that $\OR_{\mathbb{R}^d} = O(\log k)$.
	
	From Lemma \ref{rem:goedel} it follows that it is enough to provide an order on a unit cube $K = [0;1)^d$.
	For each point $x = (x_1,\dots,x_d)\in K$ we construct an infinite binary sequence $a_0a_1\dots$ as follows:
	for any $i\in (1,2,\dots,d)$ and $j\in (0,1,\dots)$ we put 
	\[
	a_{jd + i-1}= \text{the $j$-th digit in the (lex-minimal) binary representation of $x_i$.}
	\]  
	
	Let $T$  be the lexicographical order on the corresponding sequences.
	
	\begin{lemma}\label{lem:expmins}
		In $(K, T)$ there are no snakes on $2^{d+1}+1$ points with elongation greater than $8\sqrt{d}$.
	\end{lemma}
	\begin{proof}
		For any $j$ the cube $K$ is divided into $2^{jd}$ cubes of sizes
		$2^{-j}\times\dots\times2^{-j}$, we call these cubes the {\it base cubes of level $j$}.
		Note that each base cube is convex with respect to $T$.
		
		Suppose $s>2^{d+1}$ and there is a snake $x_1<_T\dots <_T x_s$ of 
		width $b$ and diameter $a > 8\sqrt{d} b$.
		Find $k$ such that $2^{-k-1} \leqslant b < 2^{-k}$.
		Points of the snake with odd indices can be covered by no more than $2^d$ base cubes of level $k$, the same holds for the points with even indices.
		
		\begin{tikzpicture}
		\draw[help lines,step = 0.5] (-0.6, -0.6) grid (3.6, 2.6);
		\draw[thick] (-0.6, -0.6) grid (3.6, 2.6);
		\draw[pattern=north west lines] (0,1) rectangle (1,2);
		\draw[pattern=north west lines] (2.5,-0.5) rectangle (3.5,0.5);
		\draw[thick, ->] (0.4,1.3) -- (3.25,0.2) -- (0.75, 1.4)--(3.35, -0.4) -- (0.75,1.75);
		\end{tikzpicture}

		By the pigeon hole principle there are two points of the snake with odd indices in one base cube $A$ of level $k$. 
		Convexity means that there must be a point with even index in the same cube $A$, but since $a - 2b > \sqrt{d} 2^{-k}$, none among the cubes of level $k$ can contain both odd points and even points.
	\end{proof}
	
	We will see later in the second claim of Lemma \ref{lem:ANfiltrationCorollary}  that the constant in Lemma \ref{lem:expmins}
	above is far from being optimal, indeed the number of the points of  snakes of large elongation, under assumption of
	Lemma \ref{lem:expmins} is at most linear in $d$.

	\section{Spaces of finite Assouad-Nagata dimension} \label{section:finitedimension}
	
	We have already mentioned that the spaces with doubling property admit orders with at most logarithmic order ratio function $\OR_{M,T}(k)$, and that among finitely generated groups only groups of polynomial growth (virtually nilpotent ones) satisfy the doubling property. 
	In this section we are going to show that any space  of finite Assouad-Nagata
	dimension admits an order with  $\OR(k) \le Const \ln k$ and provide an upper bound for the order breakpoint. We have mentioned that any group of  polynomial growth and many groups of exponential 
	growth have finite Assouad-Nagata dimension (wreath products with a base group of linear growth \cite{BrodskiyDydakLang}, Coxeter groups \cite{DranishnikovJanuszkiewicz}, relatively hyperbolic groups \cite{Hume17}, polycyclic groups \cite{HigesPeng}; see also further classes of groups mentioned in the introduction).

	Assouad-Nagata dimension can be viewed as  a (more strong, linearly controlled) version of the notion of asymptotic dimension introduced later by Gromov in \cite{gromovasymptotic}.

	\begin{definition} \label{def:nagata} Assouad-Nagata dimension of a metric space $M$ is defined as follows. 
		Let $M$ be a metric space, let  $m$  be a positive integer and $K\geqslant 1$. Suppose that for all $r>0$
		there is a family $U_r$ of subsets, such that  their union is equal to $M$ and such that the following holds
		\begin{enumerate}
			\item The diameter of any set $A \in U_r$ is at most $Kr$.
			\item Any closed (equivalently open) ball of radius  $\le r$ is contained in some set of $U_r$.
			\item Any  point of $M$ belongs to at most $m+1$ sets of $U_r$.
		\end{enumerate}
		We say in this case that Assouad-Nagata dimension of $M$ is at most $m$. To shorten the notation, we will 
		also call Assouad-Nagata dimension $AN$-dimension.
	\end{definition}
	
	If we weaken our assumption on the covering, and instead of Property $(1)$
	require that there exists some bound (not necessarily linear) for diameters of sets $A\in U_r$ in terms of $r$,
	we obtain a definition of spaces of finite {\it asymptotic dimension.}
	Such  upper bound for the diameters of $A \in U_r$ (in
	definition of $AN$-dimension, this bound is  $\le Kr$) is called
	{\it $m$-dimensional control function}.
	
	\begin{definition}[Equivalent definition of $AN$-dimension] \label{rem:def2AssouadNagata} Let $M$ be a metric space. Let us say that $M$ has $AN$-dimension at most $m$ if the following holds.
		There exists $K \geqslant 1$ such that
		for any $r$ there  exists a partition $W_r$ of the space $M$ such that all sets from $W_r$ have diameters at most $Kr$,
		and any ball of radius $r$, centered at some point $x\in M$, intersects
		at most $m+1$ sets from $W_r$.
	\end{definition}

Although the dimension $m$ does not depend on the definition, the constants $K$ in Definitions \ref{def:nagata} and \ref{rem:def2AssouadNagata} may differ.
In boths cases we call $Kr$, considered as a function of $r$, to be {\it $m$-dimensional control function}.
	
	For convenience of the reader we explain why the definitions are equivalent (with an appropriate choice $K$ for each of them).
	
	\begin{proof}

		Suppose that $M$ has $AN$-dimension
		at most $d$ with respect to the first definition.
		To explain the claim of the second definition, consider a covering $U_r$ in the first definition of Assouad-Nagata dimension. 
		For any subset $A$ from $U_r$, replace this subset by a set $A'$, which we obtain from $A$ by removing the open $r$-neighborhood of its complement.
		
		To prove that the union of $A'$ is equal to $M$ observe the following.
		For any point $x$  choose and fix a  subset $A$ from $U_r$ containing a ball of radius $r$, centered at $x$. We consider the corresponding $A'$.
		Since for any $x$ we have a set  $A'$ containing $x$, we conclude that $A'$ cover $M$.
		
		Now observe that if a ball of radius $r$ centered at $y$ intersects some $A'$, then $y$ belongs to $A$. This implies that
		any ball or radius $r$ intersects at most $m+1$ among subsets $A'$.
		Finally, observe that we can ensure that the sets are disjoint, replacing the sets by appropriate subsets.
		
		Now observe that  if $M$ admits a partition satisfying the claim of Definition \ref{rem:def2AssouadNagata}, 
		then we can  consider open $r$-neighborhoods of  the sets of the partition in this definition. 
		Let $K$ be a constant from Definition \ref{rem:def2AssouadNagata}.
		It is clear that
		the obtained sets satisfy
		the assumption of Definition \ref{def:nagata} for the same $r$ and the constant of the $m$-dimensional control function in the sense of this definition
		$K'=K+2$.
	\end{proof}
	
	As an example of a space  of finite
	Assouad-Nagata dimension recall the example
	of \cite{BrodskiyDydakLang} which shows that Assouad-Nagata dimension of 
	$\mathbb{Z}
	\wr \mathbb{Z}/2\mathbb{Z}$ is one.
	We recall that the wreath product $A \wr B$ of
	groups $A$ and $B$ is a semi-direct product of $A$ and $\sum_A B$, where $A$ acts on
	$\sum_A B$ by shifts.
	Elements of the wreath product are pairs $(a,f)$, where $a\in A$ and $f:A\to B$
	is a function with $f(x)=e_B$ for all but finitely many $x\in A$.
	We consider a standard generating set of $A\wr B$,  which
	corresponds to the union of $S_A$, and $S_B$, where $S_A \subset A \subset A\wr B$ and $S_B\subset B$ is imbedded to $\sum_A B \subset A\wr B$  by sending $B$ to the copy of $B$ in $\sum_A B$
	indexed by $e_A$. In the example below we consider the word metric on the Cayley graph with respect to this standard generating set (edges are not included).
	The argument \cite{BrodskiyDydakLang}  uses
	the fact that the kernel of the map to $\mathbb{Z}$ is zero-dimensional, and then  uses a Hurewicz type theorem for $AN$-dimension for group extensions. In the example below we describe explicitly partitions of the
	infinite  wreath product, as well as of a sequence of finite wreath products. 
	These partitions uniformly satisfy the assumption of the Definition
	\ref{rem:def2AssouadNagata} and 
	show that the spaces have
	$AN$-dimension equal to one.
	
	\begin{example}\label{ex:anwreath}[Partitions of wreath products]
		
		\begin{enumerate}
			\item
			For each $r>0$ there exists a partition $\mathcal{A}_r$
			of $G=\mathbb{Z}\wr \mathbb{Z}/2\mathbb{Z}$ such that
			diameter of every set $A
			\in \mathcal{A}$ is at
			most $6r$ and any ball of radius $r/2$ intersects at most $2$ sets.
			
			\item 
			For each $r>0$ there exists a partition $\mathcal{A}^i_r$
			of $G=\mathbb{Z}/i\mathbb{Z}\wr \mathbb{Z}/ 2\mathbb{Z}$ such that
			diameter of every set $A
			\in \mathbb{A}$ is at
			most $9r$ and any ball of radius $r/2$ intersects at most $2$ sets.
			
		\end{enumerate}
	\end{example}
	
	\begin{proof}
		(1) The proof is reminiscent of a possible argument to show that asymptotic dimension, as well as $AN$-dimension of a free group, is $1$ (see e.g. Proposition 9.8 in
		\cites{roe}).

		We partition $\mathbb{Z}$ into disjoint intervals of length $r$ and call the corresponding subset  $\{(x,f)| kr \le x < (k+1)r\}$ the {\it layer} $L_k$. It is clear that if $x\in L_k$ and $y \in L_m$ and $m-k\ge 2$ then the distance between $(x,f)$ and $(y,g)$ is $\ge r$. 
		Now we subdivide each layer $L_k$ into sets, saying that $(x,f)$ and $(y,g)$ ($kr \le  x,y <(k+1)r$) are in the same set 
		if $f(z)= g(z)$
		for any $z \notin [rk-r/2, (r+1)k+r/2]$.
		Observe that if $(x,f)$
		and $(y,g)$ are in the same layer but not in the same set, then the distance between them is $\ge r$
		since, starting from a point inside $[kr, (k+1)r]$ one needs to switch the value of the function at (at least) one point 
		of $\mathbb{Z}$ which is either $\le kr -r/2$ or $\ge (k+1)r+r/2$. 
		Observe also that diameter of any set is at most $6r$, since to go from $(x,f)$
		to any point $(y,g)$ in the same set it is sufficient to start in $x$, end at $y$, visit all points of
		$[kr -r/2,  (k+1)r+r/2]$
		(since the length of this interval is $2r$, and hence $4r$ steps suffices)  and make at most $2r$ switches.
		
		\medskip 
		\noindent (2) Similarly one can construct a  partition of a finite wreath product $\mathbb{Z}/i\mathbb{Z}\wr \mathbb{Z}/2\mathbb{Z}$. 
		If $r\ge i$, we consider the partition consisting of one set. 
		Now we assume that $r< i$. We subdivide $\mathbb{Z}/i\mathbb{Z}$ into several intervals of length $r$ and possibly one interval of length $\ge r$ and $< 2r$. 
		Each layer $L_k$ corresponding to an interval $[I_k,J_k]$ is subdivided into sets, saying that $(x,f)$ and $(y,g)$ are in the same set if $f$ and $g$ coincide on the complement of the interval $[I_k-r/2, J_k+r/2]$. As before, distinct sets 
		of the same layer are at distance at least $r$ and the sets of the layers 
		$L_k$ and $L_m$, for $k-m\ge 2$ are at distance at least $r$. The diameter of each set is at most $3r+6r$.
	\end{proof}
	
	We will return to (2) of  the example
	above in the next section,
	when we discuss the relation between  infinite $\mins$ and spectral properties of a sequence of graphs.
	In that context we will discuss
	an order, provided by
	Thm \ref{thm:nagata} for
	a disjoint union of graphs from (2) of Example \ref{ex:anwreath}.

	\begin{thm} \label{thm:nagata} If $M$ is a
		metric space of finite Assouad-Nagata dimension $m$ with $m$-dimensional control function
		at most $Kr$, then
		there exist an order $T$ on $M$ such that
		\begin{enumerate}
			\item  
			For all $k\ge 2$ it holds
			$\OR_{M,T}(k) \le C \ln k$, where a positive constant $C$ can be chosen depending on $m$ and $K$ only.
			
			\item $\mins(M,T) \leqslant 2m + 2$. Moreover, the elongations of snakes
			in $(M,T)$ on $2m+3$ points are bounded by some constant depending on $m$ and $K$ only.
		\end{enumerate}
	\end{thm}
	
	In view of Lemma \ref{rem:goedel} it is sufficient to prove the statement of 
	Theorem \ref{thm:nagata} for finite metric spaces. 
	In the lemmas below we will make a more general assumption
	that $M$ is uniformly discrete. 
	We recall that this means that 
	there exists $c>0$ such that  for all $x\ne y$ it holds $d(x,y) \ge c$.
	
	\begin{definition} [$AN$-filtrations]
		Given a metric space $M$, we say that a sequence of partitions of it $\mathcal{V}_j$, $j\in \mathbb{Z}$ is {\it $AN$-filtration} with coefficients
		$m$, $\lambda$, $D$, $\delta$,  if
		\begin{enumerate}
			\item Diameter of any set in $A \in \mathcal{V}_j$ is at most $\lambda^{j}D$. 
			\item For any $x\in M$ the ball $B(x, \lambda^{j}\delta)$ 
			can be covered by a union of  $m+1$ sets from $\mathcal{V}_j$.
			\item If $A \in \mathcal{V}_j$ and $B \in \mathcal{V}_{j'}$, $j<j'$, then either the intersection $A \cap B$ is empty, or $A \subseteq B$.
		\end{enumerate}
		In the definition above we assume that
		$m\in \mathbb{Z}_{\geq 0}$, 
		$\lambda, D, \delta \in \mathbb{R}_+$ and
		$\lambda>1$.
	\end{definition}
	
	In the definition above $j$ takes all integer values. 
	Observe that if the metric space $M$ is uniformly discrete, then for all sufficiently large $j$ the sets in $\mathcal{V}_{-j}$ are one-point sets.
	
	\begin{lemma} \label{lem:finitedimensionimpliesANfiltration}
		Let $M$ be an uniformly discrete metric space. 
		Assume that Assouad-Nagata dimension of $M$ is at most $m$ and that 
		$M$ satisfies Definition
		\ref{rem:def2AssouadNagata}
		with a linear
		constant of $m$-dimensional control  function
		at most $K$.
		Then $M$ admits an $AN$-filtration with parameters
		$m$,  $\lambda = 4K$, $\delta = \frac{1}{2}$, $D=2K$.
	\end{lemma}
	
	\begin{proof}
		Since our space is uniformly discrete, there exists $c >0$ such that for any distinct $x_1$, $x_2$ in $M$ it holds
		$d(x_1, x_2) \ge c$. 
		Observe that the constant $K$ in the definition of Assouad-Nagata dimension satisfies $K\ge 1$, and in particular $\lambda=4K \geq 4 > 1$. 
		Take the maximal  $k_0 \in \mathbb{Z}$ such that $2K \lambda^{k_0} < c$.
		For  all integers $k\le k_0$ define $\mathcal{V}_j$  to be a partition of $M$ into one-point sets.
		Observe that properties $(1), (2), (3)$ in the definition of $AN$-filtration are verified for $j,j' \le k_0$.
		
		Now suppose that we have already constructed $\mathcal{V}_j$ for all $j\le k$, $k\in \mathbb{Z}$, such that
		the properties of $AN$-filtration are verified for all $j, j'\le k$. Now we explain how to construct $\mathcal{V}_{k+1}$
		such that the properties of $AN$-filtration are verified for all $j, j' \le k+1$.
		
		Each set from $\mathcal{V}_{k+1}$ will be a union of some sets from $\mathcal{V}_{k}$.
		
		By our assumption on $M$, we know that for any $r$  there exists
		a partition $W_r$ of $M$ (into disjoint sets)  such that 
		
		1.  The diameters of the sets from $W_r$ are at most $Kr$.
		
		2.  Any ball of radius $r$, centered at some point $x\in M$, intersects
		at most $m+1$ sets from $W_r$.
		
		\noindent Consider  $r=\lambda^{k+1}$ and   the partition $W_{r}$. 
		We will group sets from $\mathcal{V}_{k}$ together using this partition.
		For any set from $W_r$ there will be one corresponding set in $\mathcal{V}_{k+1}$.
		For any set $A$ in $\mathcal{V}_k$
		we choose one of the subsets of  $W_{r}$ with a non-empty intersection with $A$, we denote this subset by $f(A)$.
		To give an idea of the size of sets $A \in \mathcal{V}_k$, we mention that by induction hypothesis the diameter of $A$ is at most $2K (4K)^k = r/2$.
		
		Now define subsets of $\mathcal{V}_{k+1}$ as the union of subsets $A$ of $\mathcal{V}_k$ with the same associated subset $f(A)$ (of $W_r$).  
		Since subsets of $\mathcal{V}_k$ are disjoint, it is clear that subsets of $\mathcal{V}_{k+1}$ are also disjoint. 
		It is also clear that for any subset $A \in \mathcal{V}_{k+1}$ and any subset 
		$B \in \mathcal{V}_j$, with $j\le k$, 
		either $B$ is contained in $A$ or they have empty intersection.
		
		Now we also observe that  there is  an upper  bound on diameters of  all subsets $A\in\mathcal{V}_{k+1}$, that guarantees that the
		Property (1) in the definition of $AN$-filtration is verified with parameters $D=2K$, $\lambda = 4K$.
		
		\[
		\sup_{A\in\mathcal{V}_{k+1}}\diam(A) \leqslant 
		\sup_{A\in W_{\lambda^{k+1}}}\diam(A) + 2\sup_{A\in\mathcal{V}_{k}}\diam(A) \leqslant K\lambda^{k+1} + 2K\lambda^{k} \leqslant 2K\lambda^{k+1}.
		\]
		
		Finally, observe that  any ball of radius  $\lambda^{k+1}$ intersects at most $m+1$ subsets of  $W_r$, 
		and hence any ball of radius 
		$\lambda^{k+1} - \sup_{A\in\mathcal{V}_{k}} \diam(A)$ intersects at most $m+1$ subsets of  $\mathcal{V}_{k+1}$.
		
		Since
		\[
		\lambda^{k+1} - \sup_{A\in\mathcal{V}_{k}}\diam(A) \geqslant 
		\lambda^{k+1} - 2K\lambda^{k} = \frac{1}{2}\lambda^{k+1},
		\]
		we can conclude that any ball of radius  $\frac{1}{2}\lambda^{k+1}$ intersects at most $m+1$ subsets 
		of $\mathcal{V}_{k+1}$. 
		Hence Property (2) on the definition of $AN$-filtration is satisfied with parameters
		$\lambda = 4K$ and $\delta =1/2$.
		
	\end{proof}
	
	\begin{lemma} \label{lem:ANfiltrationCorollary}
		Suppose that a metric space $M$ admits an $AN$-filtration $(\mathcal{V}_j)$ with coefficients $m$, $\lambda$, $D$, $\delta$, and for some order $T$ all the subsets from this $AN$-filtration are convex with respect to this order $T$. Then 
		
		\begin{enumerate}
			\item $\OR_{M,T}(k) \le C \ln k,$
			where $C$ is a constant depending on $m$, $\lambda$, $D$ and $\delta$.
			\item Elongations of snakes on $2m+3$ points are bounded by some constant $C'$, depending on $m$, $\lambda$, $D$ and $\delta$.
		\end{enumerate}
		
	\end{lemma}

	\begin{proof} 
		Consider a finite subset $X \subset M$, of cardinality $N+1$ and enumerate its points
		\[
		x_1 <_T \dots <_T x_{N+1}.
		\]
		
		Denote by $L$ the length $l_{\opt}(X)$ of a shortest path visiting all points of $X$.
		Consider some $j\in \mathbb{Z}$. Observe that the set $X$ can be partitioned into at most
		$\frac{\lambda^{-j} L}{\delta}+1$ subsets of diameter at most $\lambda^{j}\delta$.
		
		In view of Property $(2)$ in the definition of $AN$-filtration, each of these sets can be covered by union of some $m+1$ sets from $\mathcal{V}_j$. 
		We can therefore conclude that $X$ can be covered by at most 
		$(m+1)(\frac{\lambda^{-j} L}{\delta}+1)$ sets of the partition $\mathcal{V}_j$. We call these sets of $\mathcal{V}_j$ {\it parts}.
		
		\noindent Observe that if $A\in\mathcal{V}_j$, $x_i\in A$, and  $d(x_i, x_{i+1}) > \lambda^{j}D$, 
		then $x_{i+1} \notin A$.
		Since $A$ is convex with respect to $T$, all points $x_{i'}$, $i'>i$, do not belong to $A$. 
		Therefore the number of indices $i$ such that $d(x_i, x_{i+1}) > \lambda^{j}D$ is not greater than
		the total number of parts, that is, not greater than $(m+1)(\frac{\lambda^{-j} L}{\delta}+1)$.
		
		\noindent Consider $n$ and $l$ such that 
		$\lambda^{n-1} < N + 1 \leqslant  \lambda^n$, $\lambda^{l-1} < L \leqslant  \lambda^l$.
		
		Let $D_0$ be the number of indices $i$ such that 
		$\lambda^{-1}L <d(x_i, x_{i+1}) \leq L$.
		Let $D_1$ be the number of indices $i$ such that 
		$\lambda^{-2}L <d(x_i, x_{i+1}) \leq {\lambda}^{-1}L$, and so on until $D_{n-1}$.
		Finally,  we define $D_n$ as the number of indices $i$ for which $d(x_i, x_{i+1}) \leq \lambda^{-n}L$.
		
		The total length of the path with respect to our order $T$ can be estimated as 
		
		\[l_{T}(X) = \sum_{i=1}^{N}d(x_{i}, x_{i+1}) \leq L\cdot\sum_{k=0}^nD_k \lambda^{-k}.
		\]
		
		For any $k < n$, $D_k$ is not greater than the number of indices
		$i$ such that $d(x_i, x_{i+1}) > L\lambda^{-k-1}$.
		We have $L\lambda^{-k-1} > \lambda^{l - k - 2} > \lambda^{l-k-[\log_{\lambda}D]-3}D$.
		Hence, we can estimate $D_k$ as
		\[
		D_k \leq (m+1)\left(\frac{\lambda^{-(l-k-[\log_{\lambda}D]-3)}L}{\delta} + 1\right)
		\leq
		(m+1)\left(\frac{\lambda^{k+[\log_{\lambda}D]+3}}{\delta} + 1\right);
		\]
		\[
		D_k\lambda^{-k} \leq (m+1)(\lambda^{[\log_{\lambda}D+3]}/\delta + 1).
		\]
		
		\noindent Note that this number does not depend on $k$.
		We also estimate $D_n\lambda^{-n}$ as $(N+1)\lambda^{-n} \leq 1$. 
		
		\noindent Finally, 
		\[
		l_{T}(X) \leq L(n+1) (m+1)(\lambda^{[\log_{\lambda}D+3]}/\delta + 1).
		\]
		
		\noindent The first claim of the lemma follows with $C = 2(\frac{1}{\ln{\lambda}}+1)(m+1)(\lambda^{[\log_{\lambda}D+3]}/\delta +1)$.
		\medskip
		
		To prove the second claim, consider in $(M,T)$ a snake with diameter $a$ and width $b$.
		Let $k$ be such that
		$\delta\lambda^{k-1} < b \leq \delta\lambda^k$.
		Consider all points with odd indices of this snake.
		Since they form a set of
		diameter at most $\delta\lambda^k$,
		by the definition of $AN$-filtration
		this set 
		can be covered by $m+1$ subsets of $\mathcal{V}_k$. With the same
		argument we know the analogous fact
		about all even points of our snake.
		
		Since $b> \delta\lambda^{k-1}$, we know the following.
		If the elongation $a/b$ of the snake is larger than $\lambda \frac{D}{\delta} + 2$, then $a - 2b 
		>D\lambda^k$ and no two points of the snake with different parity belong to the same set from $\mathcal{V}_k$.
		In this case any of $2m+2$ covering sets 
		from $\mathcal{V}_k$ contains at most $1$ point of the snake. 
		Hence, elongations of snakes on $2m+3$ points are bounded by $\lambda(\frac{D}{\delta} + 2)$.
	\end{proof}
	
	Now we have all the ingredients to prove Theorem \ref{thm:nagata}.
	As we mentioned above, we can assume that the metric space $M$ is uniformly discrete. 
	From Lemma \ref{lem:finitedimensionimpliesANfiltration} it follows that $M$ admits an $AN$-filtration with parameters depending only on the  $m$-dimensional control function.
	Lemma \ref{le:convex} implies that for some order $T$ all the sets from this family are convex with respect to $T$, and Lemma \ref{lem:ANfiltrationCorollary} gives us needed bounds on $\OR_{M,T}$ and elongations of snakes on $2m+3$ points in $(M,T)$.
	
	\subsection{Product of two binary trees}
	In this subsection we will illustrate Theorem \ref{thm:nagata} by giving an explicit example of an order on a product of two trees.
	We encode vertices of an infinite binary tree $M$ by finite binary words $u_i\in \{0,1\}^*$, the empty word is denoted by $\varepsilon$.
	Points of the product $P$ of two binary trees $M_1$ and $M_2$
	can be coded by pairs of finite binary words $(u_1, u_2)$.
	The distance between points $(u_1,u_2)$ and $(u_3,u_4)$ in $P$ is defined as $\max(d_{M_1}(u_1, u_3), d_{M_2}(u_2, u_4))$.
	
	For given $r$ consider the following partition $W_r$ of $M$. Points $u$ and $v$ of $M$ belong to the same set  $A\in W_r$ if and only if for some non-negative integer $k$ it holds $kr \leq |u|, |v| < (k+1)r$ and the largest common prefix of $u$ and $v$ has length $\geq (k-1/2)r$. 
	Here by $|x|$ we denote the length of a word $x$, i.e. the distance to the root in $M$.
	
	This is a well-known construction of a partition  of a metric tree
	with $1$-dimensional control function $6r$, in the sense of Definition \ref{rem:def2AssouadNagata} (we referred already to this result in the proof of
	Example \ref{ex:anwreath}).
	If we take partitions $W_r$ for all integer powers of $2$, we obtain an $AN$-filtration $\mathcal{V}_k$ of $M$ with parameters 
	$\lambda = 2$, $m=1$, $D=3$ and $\delta = 1/2$.
	
	We obtain an $AN$-filtration on $P$ with parameters $\lambda = 2$, $m = 3$, $D = 3$, $\delta = 1/2$ if for each $k$ we take the partition  $\mathcal{W}_k = \{A\times B | A, B\in \mathcal{V}_k\}$. 
	(Note that $AN$-dimension of $P$ is $2$, hence from Lemma \ref{lem:finitedimensionimpliesANfiltration} there exists an $AN$-filtration of $P$ with $m=2$, but explicit description of  sets of the partitions for such filtrations would require extra work.

	Now we can describe an order $T$ on $P$ for which all sets of $\mathcal{W}_k$ are convex and consequently $\OR_{P,T}(k)$ is logarithmic in $k$ by Lemma \ref{lem:ANfiltrationCorollary}.
	
	For integers $l\geq 0$ and $k>0$  we use notation $f(l,k) = \max ([\frac{l}{2^{k+1}}]2^{k+1}-2^k, 0)$. 
	We also use notation
	$f(l,0) = l$. 
	For a point $x = (u_1, u_2)\in P$ we construct the following sequence of finite words $(v_i(x))$. 
	For any $k$, let $v_{2k}(x)$ be the prefix of $u_1$ of length $f(|u_1|, k)$ and let $v_{2k+1}(x)$ be the prefix of $u_2$ of length $f(|u_2|, k)$. Here $|u|$ denotes the length of the word $u$.
	It is clear that for $x=(u_1, u_2)$ we have $v_0(x) = u_1$, $v_1(x) = u_2$, and  $v_k(x) = \varepsilon$ for all large enough  $k$.
	Hence for any two points $x_1, x_2 \in P$, $x_1 \ne x_2$ the sequences $(v_i(x_1))$ and $(v_i(x_2))$ are distinct but coincide at all put finitely many positions.
	Finally, we put $x_1 <_T x_2$ if $v_i(x_1)<_{\rm lex} v_i(x_2)$ 
	and $i$ is the largest integer such that $v_i(x_1) \neq v_i(x_2)$.
	
	\section{Weak expansion and infinite order breakpoint}\label{sec:weakexpansion}
	
	Given a sequence of spaces $(M_\alpha, T_{\alpha})$, $\alpha \in \mathcal{A}$, we can speak about order ratio function of this sequence  by defining
	$$
	\OR(k) =  \sup_\alpha \OR_{M_\alpha, T_\alpha} (k).
	$$
	Then we can also speak about order
	breakpoint of a sequence, considering minimal $k$
	such that $\OR(k) < k$.
	
	Another way to formulate the definitions above 
	is to consider a disjoint union 
	$M = \bigcup M_\alpha$. 
	Then one can 
	define the function $\OR$
	for metric spaces with infinite distance allowed by considering only subsets with finite pairwise distances.
	In this case $\OR$ and $\mins$ of the sequence are equal
	to $\OR_M$ and $\mins(M)$.
	
	Given a finite $d$-regular graph $\Gamma$ on $n$ vertices, we consider its adjacency matrix $A$, and eigenvalues
	$$
	d=\lambda_1 \ge \lambda_2 \ge \dots \ge \lambda_n.
	$$
	It clear that $\tilde{\lambda}_i= \lambda_i/d$
	are eigenvalues
	of the normalized adjacency matrix, which is a matrix of transition probabilities of the simple random walk on $\Gamma$.
	We recall that a connected $d$-regular graph is a Ramanujan graph if $\max_{i: |\lambda_i|<d }\lambda_i  \le 2\sqrt{d-1}$.
	
	Gorodezky et al show  in \cite{gorodezkyetal}  (see also
	Theorem $2$ in Bhalgat, Chakrabarty and  Khanna  \cite{bhalgatetal})  the existence of  a sequence of bounded degree  graphs $\Gamma_i$ which have linear lower bound, in terms of cardinality of subsets,  for the competitive ratio of the universal travelling salesman problem. In our terminology, they have provided a linear lower bound for the order ratio function  of the sequence $\Gamma_i$.
	Both  \cite{gorodezkyetal}
	and  \cite{bhalgatetal} use a remarkable construction of Ramanujan graphs (constructed by Lubotzky, Philips
	and Sarnak \cite{lubotzkyetal}, see also \cite{margulis}).
	
	Now we prove that for linearity of competitive ratio in $k$, and moreover for an a priori stronger claim that
	order breakpoint is infinite,
	a milder assumption than Ramanujan condition for graphs is sufficient.
	
	Since we consider a sequence of graphs in the theorem below, we use the name of the graph as an upper index for $\lambda$, denoting $\lambda^{\Gamma_i}_j$ 
	the eigenvalue  $\lambda_j$ of $\Gamma_i$.
	
	\begin{thm}\label{thm:expanders}
		Let $\Gamma_i$ be a sequence of   finite graphs of degree $d_i\ge 3$
		on $n_i$ vertices.
		Let $T_i$ be an order on $\Gamma_i$.
		\begin{enumerate}
			
			\item  Let
			$\delta_i= (\lambda^{\Gamma_i}_1 - \lambda^{\Gamma_i}_2 )  /d_i$.
			Assume that 
			$$
			1/\delta_i  = o  \left( \frac{\log_{d_i}n_i}{\ln 
				\log_{d_i}n_i} \right).
			$$
			
			Then the order breakpoint of the sequence  $(\Gamma_i, T_i)$ is infinite. 
			
			\item Moreover, for any $k\ge 1$ any ordered graph 
			$(\Gamma_i, T_i)$ 
			admits snakes on $k$ points
			of 
			width at most 
			$$
			t_i =  \lceil\frac{24k(\ln k + \ln(2/
				\delta_i))}
			{\delta_i} \rceil
			$$
			and length  at least 
			$$
			L_i = \left[\frac{\ln n_i - 4 -2k\ln(2k) - 2k \ln (2/\delta_i)}{\ln (2d_i-1)}\right] -1
			$$
			\item In particular, for   a sequence of bounded degree expander graphs
			the following holds. If $d_i=d\ge 3$ and $\delta=\inf_i \delta_i >0$,
			then for each $k$ the graphs $(\Gamma_i, T_i)$ admit snakes on $k$ points
			of bounded width $C_k$ and of length at least $\log_{d-1}n_i -C'_k$, for some $C_k, C'_k>0$.
		\end{enumerate}
	\end{thm}
	
	\begin{proof} 
		First observe that we can easily modify our graphs to ensure that the spectrum is non-negative.
		Let $\Gamma^{+}_i$ be a graph obtained from $\Gamma_i$ by  adding $d_i$ loops in all vertices of $\Gamma_i$. Then $\Gamma^{+}_i$ is a regular graph of degree $2d_i$ (we use
		convention that each loop contributes one edge to its  vertex). 
		Observe that it is sufficient to prove that order breakpoint of a sequence  $\Gamma^{+}_i$ is infinite.
		We will also prove other claims of the theorem under assumption that the spectrum is non-negative, and then the general case will follow.
		Let $\lambda^{\Gamma^+_i}_1 = 2d_i \ge \lambda^{\Gamma^+_i}_2 \ge \dots \ge \lambda^{\Gamma^+_i}_n$ be eigenvalues of the  adjacency matrix of $\Gamma^+_i$. 
		It is clear that $\lambda^{\Gamma^+_i}_j =  \lambda_j^{\Gamma_i}+d_i \ge 0$ for all $j$. 
		In particular, $\max (\lambda_2, |\lambda_n|) =\lambda_2$. 
		Observe also that $\lambda_2^{\Gamma^+_i}/{(2d_i)} = (d_i +\lambda_2^{\Gamma_i})/{(2d_i)}$ and 
		that the spectral gap $ (
		\lambda_1^{\Gamma^+_i}-
		\lambda_2^{\Gamma^+_i} ) /{(2d_i)}
		$ is equal to  
		$\delta_i/2$. 
		
		Below we therefore assume that the graph $\Gamma_i$
		has non-negative spectrum.
		Given $i$, consider $\Gamma_i$ and choose $m_i$ starting points in $\Gamma_i$ independently with respect to the uniform distribution on the vertices of $\Gamma_i$. 
		Consider $m_i$  trajectories of independent random walks 
		with $t_i$ steps, starting in the chosen starting points. The values $m_i$ and $t_i$ will be specified later in the text.
		We will use the following observation
		\begin{claim}\label{claim:snakes}
			Let $T$ be an order on a graph $G$.
			Suppose that we have a partition of $G$
			into several convex sets with respect $T$ (we will call them intervals). Consider several $t$-step trajectories of 
			a simple random walk on $G$, and assume that the distance between starting points of any two trajectories is at least $L$. If there are at least $k$ intervals with non-empty intersection
			with two among our trajectories, then $(G,T)$ admits a snake on $k+1$ points, of width at most $t$ and length at least $L-2t$.
		\end{claim}
		
		\begin{proof}
			First observe that the distance between any two trajectories is clearly at least $L-2t$.
			
			For each among $k$ intervals choose a point $A_i$ from the first trajectory and $B_i$ from the second. Let $A_1 <_T \dots <_T A_k$ and $B_1<_T \dots <_T B_k$ be points of these trajectories. Without loss of generality, $A_1 <_T B_1$. 
			Then we can choose a snake $(A_1 <_T B_1<_T A_2<_T B_3 <_T \dots)$ on $k+1$ points that satisfies our Lemma.
			
		\end{proof}
		
		We will use Claim \ref{claim:snakes} above
		for a fixed number $k$, $G=G_i$  and  numbers $t=t_i$ and $L=L_i$, for appropriately chosen sequence $t_i$ and $L_i$.
		We will explain the proof of  (1) and (2)
		of the theorem, where for the proof of (1) we take any $k\ge 1$ and fix it for the argument below. And in (2) $k$ is the number in the formulation of this claim.
		
		We consider  the partition of the elements of $\Gamma_i$ into $N_i$ 
		sets of equal cardinality (almost equal in case $N_i$ does not divide $n_i$), convex with respect to the order $T_i$. The numbers $N_i $ will be specified later.
		In other words, we number the points of $\Gamma_i$ with respect to $T_i$:
		$g_1 \le_{T_i} g_2 \le_{T_i} \dots \le_{T_i} g_{n_i}$;
		put $r_i=[n_i/N_i]$ and consider sets
		$\Omega_1^i =\{g_1, g_2, \dots g_{r_i} \}$,
		$\Omega_2^i =\{g_{r_i+1}, g_{r_i+2}, \dots g_{2r_i}\}$, $\Omega_3^i =\{g_{2r_i+1}, g_{2r_i+2}, \dots, g_{3r_i}\}$ etc.
		
		Before we apply Claim \ref{claim:snakes} and specify
		a sequence $m_i$ of numbers of the trajectories,
		we start with a straightforward observation.
		If we choose independently $m$ points, with respect to uniform distribution on vertices of
		a graph of cardinality $n$,  of degree $d$, then the probability
		that for at least two among the chosen points
		the distance is $\le L$ is at most
		\begin{equation} \label{eq:preneravenstvo2}
		\frac{m^2d(d-1)^L}{(d-2)n}.
		\end{equation}
		
		Indeed, the number of points at distance $L$
		from a given point is at most $1+d+d(d-1)+\dots + d(d-1)^{L-1}$, and there are $m^2$ possible pairs of points chosen among $m$ points.
		
		We want to choose 
		$m_i$ and $L_i$ in such a way that
		\begin{equation} \label{eq:neravenstvo2}
		\frac{m_i^2d_i(d_i-1)^{L_i}}{(d_i-2)n_i}
		\end{equation}
		is small enough.
		In view of 
		the upper bound 
		(\ref{eq:preneravenstvo2})
		on the distance between a pair of points among $m$ points chosen at random in a graph of degree $d$,
		under  this assumption we will have a lower bound for the probability that 
		the maximum of pairwise distance between
		$m_i$ starting points of trajectories (in $\Gamma_i$)
		will be at least $L_i$.
		In order to use Claim  \ref{claim:snakes}, to find snakes of arbitrarily large elongation (and thus to prove (1) of the theorem)
		we need 
		$$
		\frac{L_i-2t_i}{t_i} \to \infty.
		$$
		We will need another condition, to estimate the number of intervals a trajectory intersects. 
		
		We put  
		\begin{equation}
		m_i = C_{N_i}^k+1.
		\end{equation}
		We use here the notation for binomial coefficients $C^k_n = \frac{n!}{k!(n-k)!}$.
		
		Since
		$m_i>C_{N_i}^k$, if each of the trajectories intersects at least $k$ convex intervals, then 
		there exist two among our $m_i$ trajectories and 
		$k$ among our $N_i$ intervals with non-empty intersection with both of these trajectories.
		
		The probability $P_i$ that one of our $m_i$ trajectories intersects less than $k$ intervals satisfies
		\begin{equation} \label{eq:corollaryThm36}
		P_i \le m_i C_{N_i}^k (\beta_i+\alpha_i)^{t_i},
		\end{equation}
		for $\beta_i =  k/N_i $ being  the density of a set which is the union of $k$  convex intervals and $\alpha_i =\lambda_{2}^{\Gamma_i}/d_i= 
		1-\delta_i$ (since we assume that the spectrum is non-negative).
		Indeed, for each among our $m_i$ trajectories
		(with starting points chosen according
		to a uniform distribution on the vertices of $\Gamma_i$)
		and fixed $k$ among
		our convex intervals, the probability to stay inside the union of these intervals is at most
		\begin{equation}\label{eq:HLWimplies}
		(\beta_i+\alpha_i)^{t_i},
		\end{equation} as follows from  \cite{HLW}[Thm 3.6]. 
		This result, that goes back to \cite{AKS},
		states the following. Let $G$ be a $d$-regular graph on $n$ vertices and 
		spectral values
		of the normalized
		adjacency matrix satisfy
		$|\tilde{\lambda}_2|, |\tilde{\lambda}_n| \le  \alpha$. Let $B$ be a subset of vertices of $G$ of cardinality
		$\le \beta n$. Then the probability that a $t$ step trajectory of a random walk,
		starting at a point chosen with respect to the uniform distribution on $G$,
		stays inside $B$ is at most $(\alpha+\beta)^t$.
		As we have already mentioned,
		for the proof of the Theorem we can assume that the spectrum of $\Gamma_i$ is non-negative, so that in this case $\max( |\lambda^{\Gamma_i}_2/d_i|, |\lambda^{\Gamma_i}_{n_i}/d_i|)=1-\delta_i$.

		We return to the proof of the theorem.
		To get a desired bound, we choose $\beta_i$ 
		of the same order as  $1-\alpha_i$, putting $\beta_i=
		k/N_i = \delta_i/2$.
		More precisely, to ensure that $N_i$ is an integer, we choose  
		\begin{equation}\label{eq:ravenstvo6}
		N_i= \lceil 2k/\delta_i \rceil.
		\end{equation}
		Then, since we can assume that $k\ge 2$ and hence $N_i \ge 4$, it holds 
		$$
		P_i   \le  (C_{N_i}^k + 1) C_{N_i}^k (\beta_i+\alpha_i)^{t_i}
		\le  2 (C_{[2k/\delta_i]}^k )^2    (1-\delta_i/2)^{t_i}
		\le 2 ((2k/\delta_i)^k)^2 \exp(-\delta_i t_i/2)=
		$$
		\begin{equation}\label{eq:preuravnenie1}
		= 2 (2k/\delta_i)^{2k} \exp(-\delta_i t_i /2),
		\end{equation}
		
		since $(1-1/x)^x \le e^{-1}$
		for any $x>1$.
		
		Take the logarithm of the previous expression.
		We want to choose $t_i$  in such a way that 
		\begin{equation}\label{eq:inequality1}
		\delta_i t_i/2   - 2k (\ln k + \ln ( 1/\delta_i)+\ln 2)-\ln2 \ge 1,    
		\end{equation}
		this guarantees
		that the probability $P_i$ that one of the trajectories intersects less than $k$ intervals is at most $\exp(-1)$.
		In particular, we can take
		\begin{equation}\label{eq:ti}
		t_i =\lceil  \frac{12k(\ln k + \ln(1/
			\delta_i))}{\delta_i} \rceil.
		\end{equation}
		
		We rewrite the upper bound in the formula (\ref{eq:preneravenstvo2})
		for the probability, that two among our $m_i$ trajectories are at distance smaller than $L_i$. 
		We take in account  the choice of $m_i$
		in the formula (\ref{eq:neravenstvo2}),
		the choice of $N_i$ in the formula
		(\ref{eq:ravenstvo6})
		and the estimate $m_i^2 \leq 2(2k/\delta_i)^{2k}$.
		We also want to assume that the obtained
		upper bound satisfies
		\begin{equation}\label{eq:qi}
		Q_i=2(2k/\delta_i)^{2k} \frac{d_i(d_i-1)^{L_i}}{(d_i-2)n_i} \le e^{-1}.
		\end{equation}
		
		Since $2e^{-1}<1$, the inequality (\ref{eq:qi})
		on $Q_i$ combined with our assumption on  $P_i$ guarantees
		that our argument,
		with  probability $\ge 1-2e^{-1}$, provides us 
		snakes of width at  most $t_i$ and of diameter at least $L_i$ in
		the ordered space
		$(\Gamma_i, T_i)$.
		
		We take the logarithm of the expression in the formula (\ref{eq:qi})
		and use that $\ln (d_i/(d_i-2)) \le \ln 3 \le 2$.
		We introduce the notation $w_i=L_i/t_i$, and
		we want therefore that $t_i$ and $w_i$
		satisfy
		\begin{equation} 
		2k\left(\ln(2k) + \ln (1/\delta_i)\right) + t_i w_i\ln (d_i-1) 
		\leq \ln n_i - 4
		\end{equation}
		
		In order to have the inequality above, to assure that $L_i$ is an integer and that we have snakes of length  $L_i$,
		we can take $w_i$ such that
		
		\begin{equation}\label{eq:wi}
		w_i  =\frac{1}{t_i}\left[\frac{\ln n_i - 4 -2k\ln(2k) - 2k \ln (1/\delta_i)}{\ln (d_i - 1)}\right].
		\end{equation}
		Thus for ordered graphs with positive spectrum we have shown existence of snake of width
		$$
		t_i =  \lceil\frac{12k(\ln k + \ln(1/
			\delta_i))}
		{\delta_i} \rceil
		$$
		and of length  at least 
		$$
		L_i = \left[\frac{\ln n_i - 4 -2k\ln(2k) - 2k \ln (1/\delta_i)}{\ln (d_i-1)}\right] -1.
		$$
		As we have mentioned  it is easy to modify the graph to ensure that the spectrum is positive.
		Since by this modification $\delta_i$ is replaced by $\delta_i/2$, $d_i$ by $2d_i$ and cardinality of the balls in $\Gamma_i$ does not change, this implies Claim (2) of the theorem for a general graph.
		
		In particular, if we have 
		a sequence of  expander graphs of fixed degree, then $1/\delta_i$ is bounded from above, hence
		there exist snakes of width at most $C$ and of length greater or equal to $\log_{d-1} n_i  -C'$, where  constants $C'$ and $C'$
		depends on $k$ and the spectral gap
		$\delta>0$, 
		$\delta= \inf_i \delta_i$.
		We have proved the third claim.
		
		To prove Claim (1) of the theorem, we want to guarantee  that the elongation of snakes tends to $\infty$, 
		and thus we need to ensure that $w_i \to \infty$, as $i \to \infty$. 
		This will follow from conditions
		\begin{equation}\label{eq:13}
		t_i \ln (d_i-1) = o (\ln n_i) 
		\end{equation}
		and
		\begin{equation}\label{eq:14}
		1 + k\ln(1/\delta_i) = o(\ln n_i).
		\end{equation}
		
		\noindent The assumption of Claim (1) of the theorem implies that $1/\delta_i = o(\ln n_i)$ and hence (\ref{eq:14}) holds.
		
		Since we can take $t_i$ as in (\ref{eq:ti}), to ensure (\ref{eq:13}) it is enough to show that 
		\begin{equation}\label{eq:dvanerva}
		\frac{\ln(1/\delta_i)}{\delta_i}  =o
		\left(
		\frac{\ln n_i}{\ln (d_i-1)}  \right). 
		\end{equation}
		
		The assumption of (1) of the Theorem implies that
		$\log_{d_i} n_i$ tends to infinity;
		and that then
		\begin{equation}
		1/\delta_i  = 
		o\left(   \frac{\log_{d_i}n_i}{\ln
			\log_{d_i}n_i} 
		\right)=
		o\left(   \frac{\log_{d_i-1}n_i}{\ln 
			\log_{d_i-1}n_i} 
		\right),
		\end{equation}  
		then 
		
		\begin{equation}
		1/\delta_i \ln(1/\delta_i)  = o \left(  \frac{\log_{d_i-1}n_i}{\ln \log_{d_i-1}n_i} \right)
		\left(\ln \log_{d_i-1}n_i - \ln \ln \log_{d_i-1}n_i \right)=
		o(\log_{d_i-1}n_i),
		\end{equation}
		
		\noindent and hence  (\ref{eq:dvanerva}) holds, and we have proved the first claim of the theorem.
	\end{proof}

	Trajectories of random walks,  for
	obtaining lower bound of the universal travelling salesman problem,  appear in \cite{gorodezkyetal} (one trajectory) and in \cite{bhalgatetal} (two trajectories). In the proof of Theorem \ref{thm:expanders}  we consider several trajectories  and an essential point of our argument is appropriate subdivision of $\Gamma_i$ into convex intervals.
	
	\begin{rem}
		Recall that a discrete version of Cheeger's inequality (see e.g. \cite{HLW}[Thm 2.4]) relates the expansion coefficient $h(\Gamma)$
		of a $d$-regular finite graph with its spectral gap:
		$$
		\frac{d-\lambda_2}{2} \le h(G) \le \sqrt{2d(d-\lambda_2)}.
		$$
		In particular, the assumption of Thm \ref{thm:expanders} is verified for a sequence of graphs of fixed degree satisfying
		$$
		\frac{1}{h(G_i)^2} = 
		o  \left( \frac{\ln n_i}{\ln
			\ln n_i} \right).
		$$
	\end{rem}

	\begin{rem}\label{rem:ANpoincare}
		Proposition $9.5$ of Hume, MacKay, Tessera \cite{humeetal} implies that a union $X$ of bounded
		degree graphs of finite $AN$-
		dimension $\Gamma_i$ and the Poincar\'{e} profile of $X$ satisfies for some $C$
		$$
		\Lambda_X^2(|\Gamma_i|) \le \frac{C |\Gamma_i|}{\log |\Gamma_i|} + C.
		$$
		Here Proposition 9.5 is applied to $X$ being a disjoint union of graphs, $\delta=1$, and 
		one uses the estimation 
		$\gamma_n(t) \le  d^{K t}$, which holds for graphs of degree $d$ and $n$-dimensional control function $Kt$.
		From Proposition  $1.2$  for $p=2$ in Bourdon \cite{bourdon}
		it follows that 
		$\Lambda^2_X(|\Gamma_i|) \geqslant |\Gamma_i|h^2(\Gamma_i)\sim |\Gamma_i|\lambda_{1,2}^{\frac{1}{2}}(\Gamma_i)$.
		In notation of Theorem \ref{thm:expanders} this means if a sequence of bounded degree
		graphs has finite $AN$-dimension, then 
		$\delta_i \le {\rm Const} (1/\ln n_i)^2$ for all $i$ and some constant ${\rm Const}$.
		In other words,
		if $1/\delta_i  = o(\ln n_i)^2$, then $AN$-dimension of this sequence is infinite.
		We are grateful to David Hume for explanations on
		\cite{humeetal} and \cite{bourdon}.
	\end{rem}

	\begin{rem}[Wreath products, lamplighter on cyclic groups]\label{rem:wreath}
		Let $\Gamma_i$ be a Cayley graph of  $\mathbb{Z}/i\mathbb{Z} \wr \mathbb{Z}/2\mathbb{Z}$, with respect to standard generators
		of this group. 
		Then $\delta_i$
		is asymptotically equivalent to  $1/i^2$, and thus is equivalent to  $1/\ln (n_i)^2$ (in fact, in this particular example all spectral
		values and the spectral measure
		for SWS  random walk are calculated, see 
		Grigorchuk, Leemann, Nagnibeda,
		Theorem  5.1 of \cite{gln} for calculation of
		spectrum and spectral measure of de Bruijn graphs and  Theorem  6.1.3 for explanation that  lamplighter graphs are particular case of de Bruijn  graphs). 
		Thus $1/\delta_i$ is equivalent to $(\ln n_i)^2$. 
		Recall that the $AN$-dimension of the disjoint union of $\mathbb{Z}/i\mathbb{Z} \wr \mathbb{Z}/2\mathbb{Z}$ is finite (see Example \ref{ex:anwreath}
		in the previous section). Thus by Theorem \ref{thm:nagata} we know that $\mins$ of this union is finite, and Claim (1) of Thm \ref{thm:expanders} does not hold.
	\end{rem}
	
	More generally, any non-trivial sequence of
	finite lamplighter graphs (and also many other sequence of wreath products) violates
	the assumption of Thm \ref{thm:expanders},
	as we discuss in the following remark.
	However, in contrast with graphs from the previous remark, many such sequences satisfy Claim (1) of the theorem. We will see it in the next section, as an application of the criterion of weak imbeddings of cubes.
	
	\begin{rem}
		There exists $C>0$ such that for any finite group $A$, of cardinality at least $2$,  the spectral gap of  $A \wr \mathbb{Z}/2 \mathbb{Z}$ satisfies
		\begin{equation}\label{eq:largerelaxtime}
		1/\delta \ge C \ln n,
		\end{equation}
		where $n$ is the cardinality of $A \wr \mathbb{Z}/2\mathbb{Z}$.
		(It is clear $n= \# A 2^{\# A})$. 
		This follows from the result of 
		Peres  and Revelle
		(\cite{peresrevelle}).
		Their result,
		under a more general assumption on the Markov chain on the wreath product (including simple random walks) states that
		$$
		t_{\rm rel} (A \wr \mathbb{Z}/2\mathbb{Z}) \ge \frac{1}{8\ln 2} t_{hit} (A),
		$$
		where {\it relaxation time} 
		$$
		t_{\rm rel} = \max_{j:\tilde{\lambda}_j <1} \frac{1}{(1-|\lambda_j/d|)}
		$$
		(in particular for random walks with non-negative spectrum
		$t_{\rm rel} = 1/(1-\lambda_2/d)$) and $t_{\rm hit}(A)$ is {\it the maximal hitting time} of a random walk on $A$, that is the maximum, over $x,y$ in the graph  of the expected time to hit $y$ starting from $x$.
		Recall that for a simple random walk on
		a finite graph $\Gamma$, hitting time
		is $\ge  \# \Gamma/2$.
		Indeed, if we fix $x\in \Gamma$, consider $T_y$ to be the time of the first visit of $y$;
		then for any trajectory $p$ starting at $x$ we have $\sum_{y \in \Gamma}T_{y} \geq \#\Gamma(\#\Gamma - 1)/2$, 
		hence $\sum_{y\in \Gamma\setminus \{x\}} \mathbb{E}_xT_y \geq \#\Gamma(\#\Gamma - 1)/2$.
		(Another result of \cite{peresrevelle} also implies that for fixed $d\ge 3$
		$A=(\Z/i \Z)^d$ there is a matching upper bound for
		$1/\delta \le {\rm Const} \ln n$).
		
		Moreover, for many wreath products  $A\wr B$ the inequality
		(\ref{eq:largerelaxtime}) holds, as can be deduced  from
		\cite{komjathyperes}[Thm 1.3]).
	\end{rem}

	\begin{rem}\label{rem:sardari}
		In claim (3) of Theorem \ref{thm:expanders} we prove
		 a logarithmic lower bound $\ln_{d-1} n_i -C_k$ for
		the length of some snakes on $k$ points in $\Gamma_i$.
		These estimates can not be significantly improved in some graphs, since this length can not be greater than the diameter:
		recall that  for any $d\ge 3$ and any $\varepsilon>0$
		the diameter of a.e. random $d$-regular graph  on $n$ vertices is at most $D$, where $D$ is the smallest integer satisfying
		$(d-1)^D \ge (2+\varepsilon) d n \ln n$
		(see \cite{BollobasdelaVega}), in particular the diameter  is close to the straightforward
		lower bound for the diameter (in terms of $n$ and
		$d$)
		$D
		=\log_{d-1} n (1+o(1))$. 
		It is known that  random 
		$d$-regular graphs form a sequence of  expander graphs, close to being Ramanujan 
		\cite{Friedman}.
		For possible Cayley graph examples with diameter close to $\log_{d-1} n$
		see
		\cite{rivinsardari}[Section 4.1], who provide a numerical evidence that the diameter of ${\rm SL}(2, \mathbb{Z}/p \mathbb{Z})$, with respect to a random generating set on $d$ generators is close to
		$\ln_{d-1} n$, as $p\to \infty$.
	\end{rem}
	
	Given a graph $\Gamma$ and an integer $l\ge 1$, denote by 
	${\rm Sc}(\Gamma, l)$ a graph obtained from $\Gamma$ by replacing each edge by a chain  of $l$ 
	edges.  We will be interested in the metric on the vertices of such graphs, and therefore to make in a regular graph we can add loops for all new vertices, so that if $\Gamma$ is a regular graph of degree $d$, then 
	the scaled graph ${\rm Sc}(\Gamma, l)$ is also regular
	of degree $d$.
	
	\begin{rem}\label{rem:nospectral}
		Let $l_i \ge 1$ be a sequence of integer numbers.
		Observe that if $\mins$
		is infinite for a sequence of graphs $\Gamma_i$, it is also infinite
		for the sequence ${\rm Sc}(\Gamma_i, l_i)$. In particular, it impossible to obtain a necessary and sufficient condition
		for $\mins=\infty$ in terms of $\delta_i$, $n_i$ and $d_i$. Indeed, consider
		a sequence of expander graphs $\Gamma_i$ and choose $l_i$ to grow rapidly. We can obtain a sequence of graphs  ${\rm Sc}(\Gamma_i, l_i)$
		with very quick decay of 
		normalized spectral gap $\delta_i$, but $\mins =\infty$. 
	\end{rem}
	
	\begin{rem}\label{rem:closetoexpandersnobs}
		On the other hand, take
		any sequence of graphs $\Gamma_i$ and choose 
		$l_i$ to grow very slow. We obtain a sequence 
		of graphs 
		$\Gamma'_i={\rm Sc}(\Gamma_i, l_i)$
		with 
		$\delta_i$ which tends to zero arbitrarily slow. It is not difficult to check that, whatever $\Gamma_i$ we take and whatever $l_i$ tending to $\infty$, the obtained graphs $\Gamma'_i$ admit orders $T_i$ (a version of ${\rm Star}$ orders)
		such that the sequence $(\Gamma'_i, T_i)$ does not contain snakes of bounded width. Thus, the assumption in (3) of the Theorem can not be weakened: no other condition on the decay of $\delta_i$ (unless $\delta_i$ is bounded away from zero) can guarantee the 
		existence of a sequence of snakes of bounded width.
	\end{rem}

	\section{Weak imbeddings of cubes and infinite order breakpoint}\label{sec:infinitegirth}
	
	In the previous section we have seen a spectral condition for a sequence of
	graphs that guarantees that order breakpoint is infinite.  In this section we prove another
	sufficient condition (for spaces or sequences of spaces) 
	in terms of weak imbeddings of cubes.
	
	The following lemma generalizes the claim of Lemma \ref{le:examplecircle} about circles. 
	We denote by $S^d$ a unit $d$-sphere in $\mathbb{R}^{d+1}$, the metrics on $S^d$ is the Euclidean metric induced from $\mathbb{R}^{d+1}$. 
	
	\begin{lemma}[Snakes in the spheres]
		\label{lem:sphere}
		Let $\varepsilon>0$ and let  $X$ be a finite  $\varepsilon$-net of the Euclidean sphere  $S^d$.  
		Let $T$ be any order on $X$. 
		Then there exist two antipodal points $x$ and $\dot{x}$ and a snake
		$(x_1<_T\dots <_T x_{d+2})$,  $x_i\in X$ for $1 \le i \le d+2$,
		such that the following holds:
		$d(x_i,x) \leq \varepsilon$ if $i$ is odd,
		$d(x_i,\dot{x}) \leq \varepsilon$ if $i$ is even.
		
		In particular, the diameter of this snake is at least  $2-2\varepsilon$ and width  $ \le 2\varepsilon$, and
		$$
		\OR_{S^d} (d+1) = d+1.
		$$
	\end{lemma}
	
	\begin{proof}
		
		Let $x$ be a point of the sphere, $s$ be a positive integer  and  $r>0$.
		Let us say that  $x$ is  {\it a tail point}  with parameters $(s,r)$, 
		if there exists a snake $x_1 <_T \dots <_T x_s$ in $X$ such that all points with odd indexes 
		are at distance at most  $r$ from $x$, and the points with even indexes are at distance at most  $r$ 
		from the point $\dot{x}$ antipodal to $x$. 
		The union of all tail points of parameters $(s,r)$ we denote by  ${\rm Tail}_{s,r}$.
		
		\begin{figure}[!htb] 
			\centering
			\includegraphics[scale=.43]{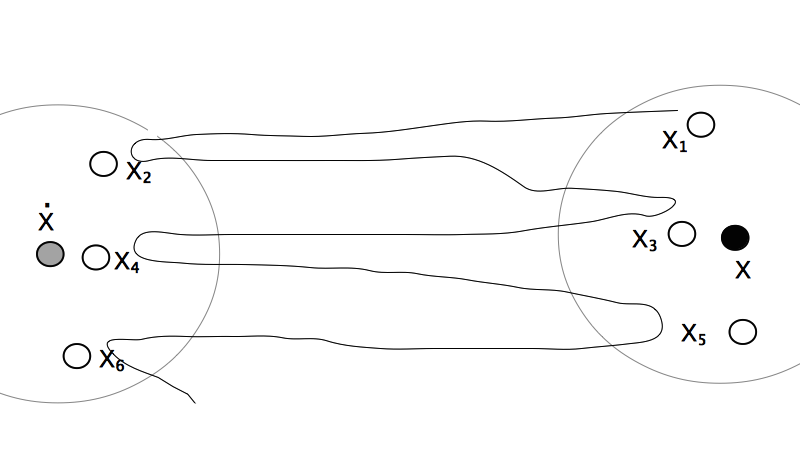}
			\caption{A tail point of a snake is presented by a black point.
			}
			\label{pic:zmejkasfera}
		\end{figure}

		\noindent Fix some positive $\delta$ much smaller then $\varepsilon$ and define open subsets $U_0,\dots,U_d$ of $S^d$ as follows.
		
		- Let $U_0$ be the set of all points at distance greater than  $\delta/2$ from ${\rm Tail}_{2, \varepsilon}$.
		
		- Let $U_1$ be the set of points at distance smaller than $\delta$ from ${\rm Tail}_{2, \varepsilon}$ but at distance greater  than $\delta/2$ from ${\rm Tail}_{3, \varepsilon + \delta}$.
		
		- Let $U_2$  be the set of points at distance smaller than $\delta$ from $\rm{Tail}_{3, \varepsilon + \delta}$,
		but at distance greater than  $\delta/2$ from $\rm{Tail}_{3, \varepsilon + 2 \delta}$, \dots
		
		- Let $U_d$ are points at distance less than $\delta$ from ${\rm Tail}_{d+1, \varepsilon + (d-1)\delta}$ and
		of distance greater than $\delta/2$ from  ${\rm Tail}_{d+2, \varepsilon + d \delta}$.
		
		If some of the sets discussed above are empty, we use the convention that a distance to an empty set is $\infty$.
		It is clear that all the subsets $U_i$ are open. 
		
		First suppose that  $U_0$, $U_1$, \dots, $U_d$  do not cover the sphere $S^d$.
		Consider a point  $x$ not belonging to their union.
		We observe that the sets  ${\rm Tail}_{2, \varepsilon}$, ${\rm Tail}_{3, \varepsilon + \delta}$, \dots,
		${\rm Tail}_{d+2, \varepsilon + d\delta}$ are at distance no more then $\delta/2$ from this point $x$. 
		In particular, the set ${\rm Tail}_{d+2, \varepsilon + d\delta}$ is not empty.
		If this happens for arbitrary small $\delta$, then by a compactness argument there exists a point $x\in S^d$ and a snake on $d+2$ points with all the odd points at distance at most $\varepsilon$ from $x$ and all the even points at distance at most $\varepsilon$ from $\dot{x}$. 
		And the claim of the lemma follows.  
		
		For the proof of the lemma it is  therefore sufficient to prove that (whatever $\delta$ we choose) the sets  $U_0$, \dots, $U_d$ can not cover the sphere.
		Suppose that it is not the case: the union of $U_i$ is equal to $S^d$. 
		We recall that by the Generalized Theorem of 
		Lusternik-Schnirelmann  (see e.g. \cite{greene}) we know that 
		\vspace{2pt}
		
		{\it	If the union of $d+1$  sets, each of the set is open or closed,  covers  $S^d$ , then for one of these sets  there exists a point $x$ of $S^d$ such that
			both $x$ and the antipodal point $\dot{x}$ belong to this set}.
		
		\vspace{2pt}
		
		We apply  this theorem for  our open sets $U_i$. We deduce that if our sets $U_0$,$U_1,$ \dots, $U_d$ cover the sphere, then there exists $k$ such 
		that $U_k$  contains a pair of antipodal points.
		
		First assume that $k=0$. Consider  $x$ such that both $x$ and $\dot{x}$ belong to $U_0$.
		Since  $X$  is  a $\varepsilon$-net of the sphere, there exists a point $a\in X$ at distance at most  $\varepsilon$ 
		from  $x$, and there exists a point $b \in X$ at distance at most $\varepsilon$ from  $\dot{x}$.
		Observe that if  $a<_Tb$, then  $x\in {\rm Tail}_{2, \varepsilon}$, and if  $a>_Tb$,  then  $\dot{x}\in {\rm Tail}_{2, \varepsilon}$. 
		We get a contradiction with the definition of the set $U_0$.
		
		Now suppose that for some $k>0$ the set $U_k$ 
		contains both  $x$ and $\dot{x}$.
		Since the distance from  $x$ to ${\rm Tail}_{k+1, \varepsilon + (k-1)\delta}$ is at most  $\delta$,
		we know that  $x \in {\rm Tail}_{k+1, \varepsilon + k\delta}$. 
		We also know that $\dot{x} \in {\rm Tail}_{k+1, \varepsilon + k\delta}$.
		Therefore, there exist snakes  $a_1 <_T \dots <_T a_{k + 1}$ and  $b_1 <_T \dots <_T b_{k+1}$  in $X$
		such that all  $a_i$ with even indexes $i$
		and all  $b_j$ with odd indexes $j$ are at distance  $< \varepsilon + k\delta$ from $\dot{x}$; and
		all  $a_j$ with odd indexes $j$ and all   $b_i$ with even indexes $i$ 
		are at distance $<\varepsilon + k\delta$ from $x$.
		
		If  $a_{k+1}<_T b_{k+1}$, then  $a_1 <_T \dots <_T a_{k + 1} <_T b_{k+1}$ is a snake and 
		we have  $x\in {\rm Tail}_{k+2, \varepsilon + k\delta}$. If  $a_{k+1} >_T b_{k+1}$, then  $\dot{x}\in {\rm Tail}_{k+2, \varepsilon + k\delta}$.
		In both cases we have obtained a contradiction with the definition of $U_k$, and we have therefore completed the proof of the lemma.
	\end{proof}

	A more combinatorial version of Lemma \ref{lem:sphere} is given in the 
	Lemma \ref{lem:triangulationsphere}
	below. We recall that the {\it octahedral triangulation} of a $d$-dimensional Euclidean sphere
	(centered at $0$)
	$S^d \subset \mathbb{R}^{d+1}$ is obtained by cutting the sphere
	by  $d+1$ pairwise distinct coordinate subspaces of dimension $d$.
	
	\begin{lemma}\label{lem:triangulationsphere}
		Consider a centrally  symmetric triangulation $K$ of the Euclidean sphere $S^d$, and assume that $K$ is a subdivision of the octahedral triangulation.
		Let $T$ be an order on vertices of $K$. 
		Then there exist two antipodal simplices $\Delta_1$ and $\Delta_2$ of this triangulation and 
		a snake on  $d+2$ points, 
		with odd indexed point being vertices of   $\Delta_1$, and even indexed points being vertices of  $\Delta_2$.
	\end{lemma}
	
	\begin{proof}
		
		Consider a triangulation $K'$, which is the barycentric subdivision of  $K$.
		By definition,  there is a one-to-one correspondence between vertices of $K'$ and simplices  of  $K$;
		two vertices of  $K'$ are adjacent  if and only if for the two corresponding simplices of $K$ one is a subset of the other.
		
		Observe that the triangulation  $K'$ is also symmetric and also is a subdivision of the  octahedral triangulation of $S^d$.
		Write in each vertex of  $K'$ an integer in the following way.
		Let $x$ and  $x'$ be two antipodal vertices of $K'$. 
		Let $\Delta$ and $\Delta'$  be corresponding (antipodal) simplices of $K$. 
		Consider maximal $s$, such that there exists a snake on  $s$ points, with all odd indexed vertices 
		in  $\Delta$ and all even indexed vertices in  $\Delta'$, or vice versa.
		Choose one of these snakes.
		If its first vertex is in  $\Delta$, assign the number  $s-1$ to the point $x$, and  $-(s-1)$ to the point $x'$.
		Analogously, if its first vertex is in  $\Delta'$, we assign $-(s-1)$ to  $x$ and $s-1$ to $x'$.
		
		Note that any two such snakes of maximal number of points can not start at opposite simplices.
		Indeed, otherwise we have snakes $x_1<_T\dots <_T x_s$ and $y_1<_T \dots <_T y_s$, 
		the points $x_s$ and $y_s$ are in opposite simplices $\Delta$ and $\Delta'$. 
		Since $s$ is the maximal length of snakes, $x_s >_T y_s$ (otherwise we can add $y_s$ to the first snake). 
		In the same way we show that $x_s <_T y_s$ and obtain a contradiction.
		
		If there is a point with assigned number $d+1$ or larger, 
		the claim of the lemma holds.
		If not, observe that   we can apply Tucker's Lemma  to a half of the sphere.
		We recall that
		this lemma (see  e.g. \cite{matucek}[Thm 2.3.1])
		claims that 
		if the vertices of a  triangulation of the $n$-ball which is antipodally
		symmetric on the boundary
		are labeled with the set $\{\pm 1,\pm 2,\dots,\pm n\}$ such that antipodal vertices on the boundary receive labels which sum to zero then some edge has labels which sum to zero. 
		Applying this lemma to a hemisphere (of $S^d$), the triangulation $K'$ and the labelling as above, we conclude that there exist two adjacent vertices of 
		$K'$ with two opposite numbers (that is, $t$ and $-t$) written in them.
		Note that this means that there are snakes on the same number of points with starting points in two antipodal simplices of $K$.
		Then at least one of them can be extended for a snake on the larger number of points. This contradiction completes the proof of
		the lemma.
	\end{proof}
	
	\begin{rem}
		A statements similar to Lemma \ref{lem:triangulationsphere}  can be obtained as a particular
		case of Zig-Zag theorem from
		Simonyi, Tardos
		\cite{zigzag}[Sec. 3.3].
		Given a topologically  $t$-chromatic graph and an ordered
		coloring of its vertices, this theorem provides a sufficient condition for existence of a full bipartite graph with colors that appear alternating
		on the two sides of the bipartite subgraph (we 
		do not recall the definition of a topologically  $t$-chromatic graph and refer to \cite{zigzag} for the exact formulation of this result).
		To study triangulations of spheres one can consider graphs with edges connecting almost antipodal
		vertices of the triangulation and obtain a version of Lemma \ref{lem:triangulationsphere}.
	\end{rem}
	
\begin{rem}
Similar constructions influence lower bounds in universal ordered protocols (also in other problems), see for example Theorem $9$ in \cite{Christodoulou1}.
\end{rem}
	
In the following definition we consider the $l_1$ metric on  cubes. 
	\begin{definition}\label{def:largescubes}
		Let us say that a metric space $(M, d_M)$
		{\it weakly contains  arbitrarily large cubes} of dimension $d$, if there exists a sequence $n_i$, tending to $\infty$ and a sequence
		of mappings of discrete cubes $f_i:[-n_i, n_i]^d\to M$ (that is, a sequence of mappings of integer points
		of 
		$[-n_i, n_i]^d$ to $M$) such that 
		$$
		\lim_{i\to \infty} \frac{{\rm minop}(f_i)}{{\rm maxne}(f_i)} = \infty.
		$$
		Here ${\rm minop}(f_i)$ is the minimal  distance $d_M(f_i(x), f_i(\dot{x}))$, where  $x, \dot{x}$ are two antipodal points on the boundary of $[-n_i,n_i]^d$, 
		and ${\rm maxne}(f_i)$
		is the maximal distance $d_M(f_i(x_1), f_i(x_2))$, where $x_1$ and $x_2$ are two neighboring points (points with integer coordinates at $l_1$-distance one) in the cube
		$[-n_i, n_i]^d$.
\end{definition}

We also say that a metric space contains {\it weakly a sequence of arbitrarily large cubes}
if for all $d\ge 1$ this space weakly contains  arbitrarily large cubes of dimension $d$.
	
	Now we explain a sufficient condition for $\OR(k) = k$.

	\begin{corol}\label{cor:cubes} 
		If a metric space $M$ weakly contains  arbitrarily large cubes of dimension $d$,
		then for any order $T$ on $M$ it holds
		$$
		OR_{M,T}(d) = d.
		$$
		
		\noindent In particular, if a metric space $M$ weakly contains a sequence of arbitrarily large cubes,
		then for any order $T$ on $M$ the order breakpoint of $(M, T)$ is infinite.
	\end{corol}
	
	\begin{proof} 
		Let $T$ be an order on $M$. 
		Assume that there exists a sequence $n_i$, tending to $\infty$, and (for each $n_i$) a mapping $f_i$ of integer points  
		of the boundary of the cube
		$[-n_i, n_i]^d$ to $M$,
		satisfying $\frac{{\rm minop}(f_i)}{{\rm maxne}(f_i)} \to \infty.
		$
		We want to  prove that $(M,T)$ admits snakes on $d+1$ points of arbitrarily large elongation.
		
		Take a sufficiently large integer $n_i$ and denote the cube $[-n_i,n_i]^d$ by $K_i$. 
		The  cube $K_i$ contains  $(2n_i+1)^d$ integer points.
		The boundary $S_i$ of $K_i$ is
		homeomorphic to a $(d-1)$-dimensional sphere, and $(2n_i+1)^d-(2n_i-1)^d$ integer points of $K_i$ belong to $S_i$. Observe that
		$S_i$ is subdivided into $2d(2n_i)^{d-1}$ unit cubes of dimension $d-1$, and this is a  subdivision of the
		(image by homeomorphism of the) octahedral partition
		of $S$.
		
		Each of the $(d-1)$-dimensional unit cubes can be divided into $(d-1)!$ simplices. 
		Hence there exists a centrally symmetric triangulation of $S_i$
		(consisting of simplices above) such that vertices of this triangulation belong to
		integer points of $K_i$ and all the simplices have diameters $d-1$.  
		
		Consider the pullback $T'$ of the order $T$  with respect to the mapping $f_i$.
		By definition $T'$ is an order on  vertices of  $S_i$ such that $x<_{T'}y \Longleftrightarrow f_i(x) <_T f_i(y)$.
		By Lemma \ref{lem:triangulationsphere} we know that there exists a pair of antipodal simplices $\Delta_1$, $\Delta_2$ and a  snake $(x_1,\dots,x_{d+1})$ on $d+1$ points oscillating between
		them.
		
		Consider the image  in $M$  (under $f_i$) of this snake. 
		Observe that $d_M(f_i(x_k),f_i(x_l)) \leqslant (d-1) {\rm maxne}(f_i) $ for any $k,l$ of the same parity.
		We assume that $i$ is large enough, so that 
		$$
		{\rm minop}(f_i)- 2 (d-1) {\rm maxne}(f_i) \ge 1.
		$$
		
		This assumption guarantees that the image of the points of the snake are distinct points in $M$. Indeed, for  the consecutive points  $f_i(x_k), f_i(x_{k+1})$ the assumption above implies
		that $d_M(f_i(x_k), f_i(x_{k+1})) \geqslant 1.$ 
		By definition of the pullback, the fact that $f_i(x_k) \ne  f_i(x_{k+1})$ and  $x_k<_{T'} x_{k+1}$, we know 
		that 
		$f_i(x_k)<_{T} f_i(x_{k+1})$. 
		We conclude therefore that
		$$
		f_i(x_1)<_{T} f_i(x_{2}) <_T \dots <_{T} f_i(x_{d})<_{T} f_i(x_{d+1}).
		$$
		
		Hence, $(f_i(x_1),\dots, (f_i(x_{d+1}))$ is indeed a snake, its width is no more than $(d-1){\rm maxne}(f_i)$ and its diameter is at least ${\rm minop}(f_i)- 2 (d-1) {\rm maxne}(f_i)$. 
		Its elongation is at least
		$$
		\frac{ {\rm minop}(f_i) - 2 (d-1) {\rm maxne}(f_i)}{{\rm maxne}(f_i)}.
		$$
		From the definition of weak imbeddings of cubes we conclude that this elongation tends to infinity, and this concludes the proof of the corollary.
	\end{proof}

	A particular case of 
	Corollary \ref{cor:cubes} is when  $G$ is such that for all $d$ there exists an uniform imbedding of $\mathbb{Z}^d$ in $G$.
	This condition holds in particular for any group $G$ that
	contains $\mathbb{Z}^\infty$ as a subgroup. 
	For example, $G = \mathbb{Z} \wr \mathbb{Z}$. Then  for any order $T$ on $G$ it holds
	$$
	\OR_{G,T}(k) = k
	$$
	for all $k$.

	We recall that some known examples of groups of infinite Assouad-Nagata dimension do not admit uniform imbeddings
	of $\mathbb{Z}^d$. For example, if $G=\mathbb{Z}^2\wr A$, $A$ is a finite group  of cardinality $\ge 2$, then
	Assouad-Nagata dimension of $G$ is infinite,
	see \cite{nowak}, who studied also other amenable wreath products; for a general case of any base group of superlinear growth see \cite{BrodskiyDydakLang}[Corr 5.2]. On the other hand, the asymptotic dimension of $G$ is $2$, see \cite{nowak}, see also   \cite{BrodskiyDydakLang}[Thm 4.5]
	for upper bounds on dimensional control function. In particular $G$ can not contain uniform images of $\mathbb{Z}^3$.
	Wreath product examples mentioned above  are discussed in the following subsection.
	
	\subsection{Wreath products}

	\begin{lemma}\label{lem:wcubeswr}
		[Sequences of arbitrarily large cubes
		in wreath products]
		
		\noindent Let $G =A \wr B$, where $A$ is a finitely generated group of super-linear growth and 
		$B$ is a finitely generated group of cardinality at least two, then $G$ contains weakly a sequence of arbitrarily large cubes.
	\end{lemma}
	
	The argument we explain below works for all $B$, but we want
	to point out that our main interest when $B$ is finite.
	(If $B$ is infinite and finitely generated, observe that for all $d\ge 1$ $G$ contains uniformly $\mathbb{Z}_+^d$, and the claim of the lemma follows).
	The proof below is
	reminiscent of
	an argument  of  Theorem 4.1 in  \cite{BrodskiyDydakLang}.
	
	\begin{proof}
		Observe that the claim of the lemma does not
		depend on the generating set in the wreath product.
		Fix some generating set $S_A$ of $A$, a set $S_B$ of $B$ and consider a standard generating set
		$S$ of $G$, which one can identify with the union of $S_A$ and $S_B$.
		Fix $d\ge 1$. We are going to prove that $G$ contains weakly arbitrarily large cubes of dimension $d$.
		
		We first observe that for infinitely many $n$ there exist $d$ disjoint
		subsets $\Omega^n_1$, $\Omega^n_2$,  \dots, $\Omega^n_d$ inside the ball $B_{A,S_A}(e,n)$ of radius $n$ such that the following holds.
		The cardinality of each set $\Omega_i$  is $n$ and for each $i$, $1 \le i \le d$,
		the length of
		any path visiting all  points of $\Omega^n_i$ is  $\ge \frac{1}{6} (n-1) \sqrt{n}$.
		
		We have already mentioned an elementary case of Polynomial Growth theorem of
		\cite{justin71}. The result of this paper states 
		that if $A$ is not virtually cyclic, then the growth function of $A$ satisfies
		$v_{A,S_A}(n) \ge n(n+1)/2$.
		In this case  there exist
		infinitely many $n$ such that 
		$v_{A,S_A}(n+\sqrt n)/v_{A,S_A}(\sqrt{n}/6) \ge n$.
		Take such $n$ and put $\varepsilon = \sqrt{n}/2$, choose an $(\varepsilon, \varepsilon)$-net of the ball $B_{A,S_A}(e,n)$ in $A$. 
		Place a ball of radius $\varepsilon/3$ around each point of the net.
		In each ball choose $d$ distinct elements. 
		If $n$ is sufficiently large, then $v_{A,S_A}(\sqrt{n}/6)$ is greater
		than $d$, and therefore  such choice is possible.  
		Observe also that the number of the points
		in such net is at least $v_{A,S_A}(n+\sqrt n)/v_{A,S_A}(\sqrt{n}/6)\geq n$. 
		
		In each ball choose one point and denote the 
		union of
		the chosen points $\Omega^n_1$. Choose one more 
		point in each ball and denote the union of these points by $\Omega^n_2$ and so on.
		It is clear that the distance between any two points of $\Omega_i^n$ is at least $\sqrt{n}/6$, and hence 
		the length of any path visiting all  points of $\Omega^n_i$ is at least $\frac{1}{6} (n-1) \sqrt{n}$.
		
		\begin{figure}[!htb] 
			\centering
			\includegraphics[scale=.85]{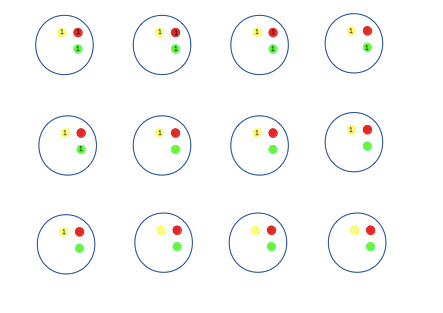}
			\caption{The sets $\Omega^n_i$, here shown for  $n=12$, $d=3$, $i=1, 2, 3$. The points of the same colour correspond to the same set.
				The image of $(9,3,5)$ under $\rho$ is shown, this is a configuration taking value $1$ in the 
				first $9$ yellow points, in the first $3$ red points and in the first
				$5$ green points.}
			\label{pic:lc}
		\end{figure}
		
		Enumerate the points of $(\varepsilon, \varepsilon)$-net in an arbitrary order and consider the restriction of this order to $\Omega_i^n$.
		Fix a non-identity element $b$ in the already fixed generating set of   $B$.
		Consider  a map  $\rho: [0,1,\dots,n]^d \to A \wr B$.
		A point with coordinates $(z_1, z_2, \dots, z_d)$ is sent to $(e_A,f)$ where $f$ is the configuration  which takes value
		$b$ in the first $z_i$ elements of $\Omega_i$, $1\le i \le d$, and takes value $e_B$ elsewhere.
		
		See Picture \ref{pic:lc}. The picture shows a possible
		choice for $A=\mathbb{Z}^2$,
		$B=\mathbb{Z}/2\mathbb{Z}$. 
		We use an additive notation
		in $B$, and a non-identity element $b$ is denoted by $1$.
		
		If  $u, v$ are points in the cube
		at distance $1$ in  $l_1$-metric of $\mathbb{Z}^d$, then the distance between $\rho(u)$ and $\rho(v)$ in the word metric of the wreath product is at most $2n+1$. 
		Indeed, observe that in this case the values of the configuration
		of $a$ and that of $b$ differ at one point, which we denote by $x$. To make this change, it is sufficient to go from identity to $x$, to make a switch at go back to the identity.
		
		Observe that if we take two antipodal points $w, \dot{w}$ in the boundary of the cube, then the distance between $\rho(w)$ and $\rho(\dot{w})$ is at least $\frac{1}{6} \sqrt{n} (n-1)$.
		
		Indeed, observe that if we have a pair of antipodal points on the boundary of the cube, then there exists $i$ such that the $i$-th coordinate of one of them is $n$ and of the other is $0$.
		Therefore, 
		To move from between the images of these antipodal points, 
		we need the to visit all points
		of $\Omega^n_i$, and hence the length of such path is at least $\sqrt{n} (n-1)/6$.
		
		This completes the proof of the lemma.
	\end{proof}

	\begin{corol}[Order ratio function for wreath products]\label{cor:wreathgap}
		Let $G=A \wr B$ be a wreath product of $A$ and $B$,  where $\#A = \infty, \#B > 1$, $A$ and $B$ are finitely generated.
		Then either we have a 
		logarithmic upper bound for the order
		ratio function (and this happens if and only if  $G$ has
		a finite $AN$-dimension), or the order ratio function is linear (moreover, 
		the order breakpoint  is infinite).
	\end{corol}
	\begin{proof}
		Indeed, by a result of \cite{BrodskiyDydakLang}[Thm 5.1 and Cor 5.2]
		$AN$-dimension of $A\wr B$ is finite if and only if $A$ is of linear growth and $B$ is finite.
		The cited theorem deals with the case when $B$ is finite, and it is straightforward that the $AN$-dimension is infinite when $A$
		and $B$ are both infinite, and the dimension is finite when $A$ and $B$
		are both finite.
		
		By Thm \ref{thm:nagata}
		we know that if $AN$-dimension is finite, then there is a logarithmic bound for the order ratio function. On the other hand, if $A$ is of super-linear growth or if $A$ and $B$ are infinite, 
		we know by Lemma \ref{lem:wcubeswr} that $A\wr B$ contains weakly arbitrarily large cubes (of arbitrarily large dimension), and hence by Corollary \ref{cor:cubes} the order breakpoint is infinite.
		
	\end{proof}

	\subsection{Product of tripods}
	
	In Thm \ref{thm:nagata} we proved that if a metric space $M$ has $AN$-dimension $d$ then $\mins(M) \leq 2d+2$.
	Now we will show that this estimate is close to optimal.
	
	\begin{prop}
		Let $M$ be a Cartesian product of $d$ tripods.
		Then for any order $T$ on $M$ it holds $\OR_{M,T}(2d) = 2d$;
		in other words, $\mins(M) \geq 2d+1$.
	\end{prop}
	
	\begin{proof}
		There exists a continuous map $f:S^{2d-1} \to M$ such that any two opposite points in $S^{2d-1}$ map to different points.
		
		Indeed, $S^{2d-1}$ is homeomorphic to the boundary of the product $D^d$ of $d$ $2$-dimensional unit disks.
		There is a continuous mapping $h$ from a unit disk to tripod such that $h$ maps any two antipodal points on the boundary of $D$ to points at distance $1$.

		\begin{tikzpicture}[scale = 0.5, every node/.style={scale=0.8}]
		\draw[pattern=north west lines] (0,1) circle(3);
		\foreach \x in {0,...,5}{
			\node at ({3.2*cos(60 * \x + 30)}, {3.2*sin(60 * \x + 30)+1}) {$A_{\x}$};
		}
		
		\draw[->][thick] (6.5,1.1) -- (8.5,1.1);
		
		\begin{scope}[shift = {(13,0)}]
		\draw[ultra thick] (0,0) -- (0,3);
		\draw[ultra thick] (0,0) -- (2.5,-1);
		\draw[ultra thick] (0,0) -- (-2.5,-1);
		\draw (-0.2,3.2) node[above]{$A_1$} -- (-0.2, 0.2) node[above left]{$A_2$} --
		(-2.6, -0.8) to[out = -130, in = 190] (-2.4, -1.2) node[below left] {$A_3$} -- 
		(0, -0.25)node[below] {$A_4$} -- (2.4, -1.2) node[below right] {$A_5$}  to [out = -20, in = -10] (2.6, -0.8)
		--(0.2, 0.2) node[above right] {$A_0$} -- (0.2, 3.2) to [out = 90, in = 90] cycle;
		\end{scope}
		
		\end{tikzpicture}
		
		Denote the mapping $h^d:D^d \to M$ by $f$.
		If $x$ and $x'$ are two opposite points in the boundary of $D^d$ then their projections to one of the disks are opposite points of boundary of this disk, and projections of $f(x)$ and $f(x')$ to the corresponding tripod are two points at distance 1.
		
		Let $T$ be an order on $M$ and let $T'$ be a pullback of $T$ to $S^{2d-1}$.
		By Lemma \ref{lem:sphere} for any $\varepsilon$ there are two opposite points $x$ and $x'$ of $S^{2d-1}$ and a snake $x_1 <_{T'} \dots <_{T'} x_{2d+1}$ such that points with odd indices are in $B(x,\varepsilon)$ and points with even indices are in $B(x',\varepsilon)$. 
		It is clear that points $f(x_i)$ form a snake of large elongation in $(M,T)$, 
		so $\OR_{M,T}(2d) = 2d$. 
	\end{proof}


\begin{thebibliography}{}
		
		\bib{AKS}{article}{
			author={Ajtai, M.},
			author={Koml\'{o}s, J.},
			author={Szemer\'{e}di, E.},
			title={Sorting in $c\,{\rm log}\,n$ parallel steps},
			journal={Combinatorica},
			volume={3},
			date={1983},
			number={1},
			pages={1--19},
			issn={0209-9683},
			review={\MR{716418}},
			doi={10.1007/BF02579338},
		}
		
		
		\bib{Assouad82}{article}{
			author={Assouad, Patrice},
			title={Sur la distance de Nagata},
			language={French, with English summary},
			journal={C.~R.~Acad. Sci. Paris S\'{e}r. I Math.},
			volume={294},
			date={1982},
			number={1},
			pages={31--34},
			issn={0249-6321},
			review={\MR{651069}},
		}

		
		\bib{bartholdiplatzman82}{article}{
			author={Bartholdi, J. J., III},
			author={Platzman, L. K.},
			title={An $O(N\,{\rm log}\,N)$ planar travelling salesman heuristic based
				on spacefilling curves},
			journal={Oper. Res. Lett.},
			volume={1},
			date={1982},
			number={4},
			pages={121--125},
			issn={0167-6377},
			review={\MR{687352}},
			doi={10.1016/0167-6377(82)90012-8},
		}
		
		\bib{bartholdiplatzman89}{article}{
			author={Platzman, Loren K.},
			author={Bartholdi, John J., III},
			title={Spacefilling curves and the planar travelling salesman problem},
			journal={J. Assoc. Comput. Mach.},
			volume={36},
			date={1989},
			number={4},
			pages={719--737},
			issn={0004-5411},
			review={\MR{1072243}},
			doi={10.1145/76359.76361},
		}
		
		\bib{bertsimasgrigni}{article}{
			author={Bertsimas, Dimitris},
			author={Grigni, Michelangelo},
			title={Worst-case examples for the spacefilling curve heuristic for the
				Euclidean traveling salesman problem},
			journal={Oper. Res. Lett.},
			volume={8},
			date={1989},
			number={5},
			pages={241--244},
			issn={0167-6377},
			review={\MR{1022839}},
			doi={10.1016/0167-6377(89)90047-3},
		}
		
		\bib{bhalgatetal}{article}{
			author={Bhalgat, Anand},
			author={Chakrabarty, Deeparnab},
			author={Khanna, Sanjeev},
			title={Optimal lower bounds for universal and differentially private
				Steiner trees and TSPs},
			conference={
				title={Approximation, randomization, and combinatorial optimization},
			},
			book={
				series={Lecture Notes in Comput. Sci.},
				volume={6845},
				publisher={Springer, Heidelberg},
			},
			date={2011},
			pages={75--86},
			review={\MR{2863250}},
		}
		
		
		\bib{BollobasdelaVega}{article}{
			author={Bollob\'{a}s, B.},
			author={Fernandez de la Vega, W.},
			title={The diameter of random regular graphs},
			journal={Combinatorica},
			volume={2},
			date={1982},
			number={2},
			pages={125--134},
			issn={0209-9683},
			review={\MR{685038}},
			doi={10.1007/BF02579310},
		}

		
		\bib{bourdon}{article}{
			author={Bourdon, Marc},
			title={Un th\'{e}or\`eme de point fixe sur les espaces $L^p$},
			language={French, with English summary},
			journal={Publ. Mat.},
			volume={56},
			date={2012},
			number={2},
			pages={375--392},
			issn={0214-1493},
			review={\MR{2978328}},
			doi={10.5565/PUBLMAT\_56212\_05},
		}

		
		\bib{BrodskiyDydakLang}{article}{
			author={Brodskiy, N.},
			author={Dydak, J.},
			author={Lang, U.},
			title={Assouad-Nagata dimension of wreath products of groups},
			journal={Canad. Math. Bull.},
			volume={57},
			date={2014},
			number={2},
			pages={245--253},
			issn={0008-4395},
			review={\MR{3194169}},
			doi={10.4153/CMB-2013-024-8},
		}

		
		\bib{BuyaloDranishnikovSchroeder}{article}{
			author={Buyalo, Sergei},
			author={Dranishnikov, Alexander},
			author={Schroeder, Viktor},
			title={Embedding of hyperbolic groups into products of binary trees},
			journal={Invent. Math.},
			volume={169},
			date={2007},
			number={1},
			pages={153--192},
			issn={0020-9910},
			review={\MR{2308852}},
			doi={10.1007/s00222-007-0045-2},
		}
		
		
		\bib{Christodoulou1}{article}{
			author={Christodoulou, George},
			author={Sgouritsa, Alkmini},
			title={Designing networks with good equilibria under uncertainty},
			conference={
				title={Proceedings of the Twenty-Seventh Annual ACM-SIAM Symposium on
					Discrete Algorithms},
			},
			book={
				publisher={SODA  2016,  Arlington,  VA, USA},
			},
			date={2016},
			pages={72-89},
		}
		
		
		\bib{Christodoulou}{article}{
			author={Christodoulou, George},
			author={Sgouritsa, Alkmini},
			title={An improved upper bound for the universal TSP on the grid},
			conference={
				title={Proceedings of the Twenty-Eighth Annual ACM-SIAM Symposium on
					Discrete Algorithms},
			},
			book={
				publisher={SIAM, Philadelphia, PA},
			},
			date={2017},
			pages={1006--1021},
			review={\MR{3627793}},
			doi={10.1137/1.9781611974782.64},
		}
	
		
		\bib{DranishnikovJanuszkiewicz}{article}{
			author={Dranishnikov, A.},
			author={Januszkiewicz, T.},
			title={Every Coxeter group acts amenably on a compact space},
			booktitle={Proceedings of the 1999 Topology and Dynamics Conference (Salt
				Lake City, UT)},
			journal={Topology Proc.},
			volume={24},
			date={1999},
			number={Spring},
			pages={135--141},
			issn={0146-4124},
			review={\MR{1802681}},
		}
		
	
		
		\bib{ErschlerMitrofanov1}{article}{
			author={Erschler, Anna},
			author={Mitrofanov, Ivan},
			title={Spaces that can be ordered effectively: virtually free groups and hyperbolicity},
			year={2021},
			status={preprint, https://arxiv.org/abs/2011.01732 }
		}	
		
		
		\bib{Friedman}{article}{
			author={Friedman, Joel},
			title={A proof of Alon's second eigenvalue conjecture and related
				problems},
			journal={Mem. Amer. Math. Soc.},
			volume={195},
			date={2008},
			number={910},
			pages={viii+100},
			issn={0065-9266},
			isbn={978-0-8218-4280-5},
			review={\MR{2437174}},
			doi={10.1090/memo/0910},
		}
		
		
		\bib{eades}{article}{
			author={Eades, Patrick Fintan},
			title={Uncertainty Models in Computational Geometry},
			date={2020},
			journal={Phd thesis, University of Sydney, Retrieved from https://ses.library.usyd.edu.au/bitstream/handle/2123/23909/thesis\%20revised\%20(1).pdf },
		}
		
		\bib{eades2}{article}{
			author = {Eades, Patrick},
			author= {Mestre, Juli\'{a}n},
			title = {An optimal lower bound for hierarchical universal solutions for TSP on the plane},
			journal = {Computing and combinatorics. 26th international conference, COCOON 2020, Atlanta, GA, USA, August 29--31, 2020. Proceedings},
			pages = {222--233},
			date = {2020},
			publisher = {Cham: Springer},
		}
		
		\bib{gln}{article}{
			author={Grigorchuk, R.},
			author={Leemann, P.-H.},
			author={Nagnibeda, T.},
			title={Lamplighter groups, de Brujin graphs, spider-web graphs and their
				spectra},
			journal={J. Phys. A},
			volume={49},
			date={2016},
			number={20},
			pages={205004, 35},
			issn={1751-8113},
			review={\MR{3499181}},
			doi={10.1088/1751-8113/49/20/205004},
		}

		
		\bib{gorodezkyetal}{article}{
			author={Gorodezky, Igor},
			author={Kleinberg, Robert D.},
			author={Shmoys, David B.},
			author={Spencer, Gwen},
			title={Improved lower bounds for the universal and {\it a prior}i TSP},
			conference={
				title={Approximation, randomization, and combinatorial optimization},
			},
			book={
				series={Lecture Notes in Comput. Sci.},
				volume={6302},
				publisher={Springer, Berlin},
			},
			date={2010},
			pages={178--191},
			review={\MR{2755835}},
		}

		
		
		\bib{greene}{article}{
			author={Greene, Joshua E.},
			title={A new short proof of Kneser's conjecture},
			journal={Amer. Math. Monthly},
			volume={109},
			date={2002},
			number={10},
			pages={918--920},
			issn={0002-9890},
			review={\MR{1941810}},
			doi={10.2307/3072460},
		}
		
		
		\bib{gromovpolynomial}{article}{
			author={Gromov, Mikhael},
			title={Groups of polynomial growth and expanding maps},
			journal={Inst. Hautes \'{E}tudes Sci. Publ. Math.},
			number={53},
			date={1981},
			pages={53--73},
			issn={0073-8301},
			review={\MR{623534}},
		}

		
		\bib{gromovasymptotic}{article}{
			author={Gromov, M.},
			title={Asymptotic invariants of infinite groups},
			conference={
				title={Geometric group theory, Vol. 2},
				address={Sussex},
				date={1991},
			},
			book={
				series={London Math. Soc. Lecture Note Ser.},
				volume={182},
				publisher={Cambridge Univ. Press, Cambridge},
			},
			date={1993},
			pages={1--295},
			review={\MR{1253544}},
		}
		

		
		\bib{HaglundWise08}{article}{
			author={Haglund, Fr\'{e}d\'{e}ric},
			author={Wise, Daniel T.},
			title={Special cube complexes},
			journal={Geom. Funct. Anal.},
			volume={17},
			date={2008},
			number={5},
			pages={1551--1620},
			issn={1016-443X},
			review={\MR{2377497}},
			doi={10.1007/s00039-007-0629-4},
		}
		

		\bib{higes}{article}{
			author={Higes, J.},
			title={Assouad-Nagata dimension of locally finite groups and asymptotic
				cones},
			journal={Topology Appl.},
			volume={157},
			date={2010},
			number={17},
			pages={2635--2645},
			issn={0166-8641},
			review={\MR{2725356}},
			doi={10.1016/j.topol.2010.07.015},
		}
		
		
		\bib{HigesPeng}{article}{
			author={Higes, J.},
			author={Peng, I.},
			title={Assouad-Nagata dimension of connected Lie groups},
			journal={Math. Z.},
			volume={273},
			date={2013},
			number={1-2},
			pages={283--302},
			issn={0025-5874},
			review={\MR{3010160}},
			doi={10.1007/s00209-012-1004-1},
		}
		
		\bib{kapovich}{article}{
			author={Kapovich, Michael},
			title={Triangle inequalities in path metric spaces},
			journal={Geom. Topol.},
			volume={11},
			date={2007},
			pages={1653--1680},
			issn={1465-3060},
			review={\MR{2350463}},
			doi={10.2140/gt.2007.11.1653},
		}
		
		\bib{matucek}{book}{
			author={Matou\v{s}ek, Ji\v{r}\'{\i}},
			title={Using the Borsuk-Ulam theorem},
			series={Universitext},
			note={Lectures on topological methods in combinatorics and geometry;
				Written in cooperation with Anders Bj\"{o}rner and G\"{u}nter M.~Ziegler},
			publisher={Springer-Verlag, Berlin},
			date={2003},
			pages={xii+196},
			isbn={3-540-00362-2},
			review={\MR{1988723}},
		}
		

		
		\bib{HLW}{article}{
			author={Hoory, Shlomo},
			author={Linial, Nathan},
			author={Wigderson, Avi},
			title={Expander graphs and their applications},
			journal={Bull. Amer. Math. Soc. (N.S.)},
			volume={43},
			date={2006},
			number={4},
			pages={439--561},
			issn={0273-0979},
			review={\MR{2247919}},
			doi={10.1090/S0273-0979-06-01126-8},
		}
		
		
		
		\bib{Hume17}{article}{
			author={Hume, David},
			title={Embedding mapping class groups into a finite product of trees},
			journal={Groups Geom. Dyn.},
			volume={11},
			date={2017},
			number={2},
			pages={613--647},
			issn={1661-7207},
			review={\MR{3668054}},
			doi={10.4171/GGD/410},
		}
		
		\bib{humeetal}{article}{
			author={Hume, David},
			author={Mackay, John},
			author={Tessera, Romain},
			title={Poincar\'{e} profiles of groups and spaces},
			journal={Rev. Mat. Iberoam.},
			volume={36},
			date={2020},
			number={6},
			pages={1835--1886},
			issn={0213-2230},
			review={\MR{4163985}},
			doi={10.4171/rmi/1184},
		}
		
		\bib{jiadoubling}{article}{
			author={Jia, Lujun},
			author={Lin, Guolong},
			author={Noubir, Guevara},
			author={Rajaraman, Rajmohan},
			author={Sundaram, Ravi},
			title={Universal approximations for TSP, Steiner tree, and set cover},
			conference={
				title={STOC'05: Proceedings of the 37th Annual ACM Symposium on Theory
					of Computing},
			},
			book={
				publisher={ACM, New York},
			},
			date={2005},
			pages={386--395},
			review={\MR{2181640}},
			doi={10.1145/1060590.1060649},
		}
		
		
		\bib{justin71}{article}{
			author={Justin, Jacques},
			title={Groupes et semi-groupes \`a croissance lin\'{e}aire},
			language={French},
			journal={C.~R.~Acad. Sci. Paris S\'{e}r. A-B},
			volume={273},
			date={1971},
			pages={A212--A214},
			issn={0151-0509},
			review={\MR{289689}},
		}
		
		\bib{komjathyperes}{article}{
			author={Komj\'{a}thy, J\'{u}lia},
			author={Peres, Yuval},
			title={Mixing and relaxation time for random walk on wreath product
				graphs},
			journal={Electron. J. Probab.},
			volume={18},
			date={2013},
			pages={no. 71, 23},
			review={\MR{3091717}},
			doi={10.1214/EJP.v18-2321},
		}
		

		
		\bib{LangSchlichenmaier}{article}{
			author={Lang, Urs},
			author={Schlichenmaier, Thilo},
			title={Nagata dimension, quasisymmetric embeddings, and Lipschitz
				extensions},
			journal={Int. Math. Res. Not.},
			date={2005},
			number={58},
			pages={3625--3655},
			issn={1073-7928},
			review={\MR{2200122}},
			doi={10.1155/IMRN.2005.3625},
		}
		
		
		
		\bib{lubotzkyetal}{article}{
			author={Lubotzky, A.},
			author={Phillips, R.},
			author={Sarnak, P.},
			title={Ramanujan graphs},
			journal={Combinatorica},
			volume={8},
			date={1988},
			number={3},
			pages={261--277},
			issn={0209-9683},
			review={\MR{963118}},
			doi={10.1007/BF02126799},
		}
		
		\bib{MannHowgroupsgrow}{book}{
			author={Mann, Avinoam},
			title={How groups grow},
			series={London Mathematical Society Lecture Note Series},
			volume={395},
			publisher={Cambridge University Press, Cambridge},
			date={2012},
			pages={x+199},
			isbn={978-1-107-65750-2},
			review={\MR{2894945}},
		}
		
		\bib{margulis}{article}{
			author={Margulis, G. A.},
			title={Explicit group-theoretic constructions of combinatorial schemes
				and their applications in the construction of expanders and
				concentrators},
			language={Russian},
			journal={Problemy Peredachi Informatsii},
			volume={24},
			date={1988},
			number={1},
			pages={51--60},
			issn={0555-2923},
			translation={
				journal={Problems Inform. Transmission},
				volume={24},
				date={1988},
				number={1},
				pages={39--46},
				issn={0032-9460},
			},
			review={\MR{939574}},
		}
		
		
		\bib{Nagata58}{article}{
			author={Nagata, J.},
			title={Note on dimension theory for metric spaces},
			journal={Fund. Math.},
			volume={45},
			date={1958},
			pages={143--181},
			issn={0016-2736},
			review={\MR{105081}},
			doi={10.4064/fm-45-1-143-181},
		}
		
		\bib{nowak}{article}{
			author={Nowak, Piotr W.},
			title={On exactness and isoperimetric profiles of discrete groups},
			journal={J. Funct. Anal.},
			volume={243},
			date={2007},
			number={1},
			pages={323--344},
			issn={0022-1236},
			review={\MR{2291440}},
			doi={10.1016/j.jfa.2006.10.011},
		}
		
\bib{osajda}{article}{
   author={Osajda, Damian},
   title={Small cancellation labellings of some infinite graphs and
   applications},
   journal={Acta Math.},
   volume={225},
   date={2020},
   number={1},
   pages={159--191},
   issn={0001-5962},
   review={\MR{4176066}},
   doi={10.4310/acta.2020.v225.n1.a3},
}
		
		\bib{peresrevelle}{article}{ author={Peres, Yuval}, author={Revelle, David}, title={Mixing times for random walks on finite lamplighter groups}, journal={Electron. J. Probab.}, volume={9}, date={2004}, pages={no. 26, 825--845}, issn={1083-6489}, review={\MR{2110019}}, doi={10.1214/EJP.v9-198},}
		
		\bib{rivinsardari}{article}{
			author={Rivin, Igor},
			author={Sardari, Naser T.},
			title={Quantum chaos on random Cayley graphs of ${\rm SL}_2[\mathbb{ Z}/p\mathbb
				{Z}]$},
			journal={Exp. Math.},
			volume={28},
			date={2019},
			number={3},
			pages={328--341},
			issn={1058-6458},
			review={\MR{3985838}},
			doi={10.1080/10586458.2017.1403982},
		}
		
		
		\bib{roe}{book}{
			author={Roe, John},
			title={Lectures on coarse geometry},
			series={University Lecture Series},
			volume={31},
			publisher={American Mathematical Society, Providence, RI},
			date={2003},
			pages={viii+175},
			isbn={0-8218-3332-4},
			review={\MR{2007488}},
			doi={10.1090/ulect/031},
		}
		

		
\bib{Sledd1}{article}{
   author={Sledd, Levi},
   title={Assouad-Nagata dimension of finitely generated $C'(1/6)$ groups},
   journal={Topology Appl.},
   volume={299},
   date={2021},
   pages={Paper No. 107634, 17},
   issn={0166-8641},
   review={\MR{4265774}},
   doi={10.1016/j.topol.2021.107634},
}
		
		
		\bib{Sledd}{article}{
			author={Sledd, Levi},
			title={Asymptotic and Assouad-Nagata dimension of finitely generated groups and their subgroups},
			year={2020},
			status={preprint, https://arxiv.org/abs/2010.03709}
		}
		
		

		
		\bib{zigzag}{article}{
			author={Simonyi, G\'{a}bor},
			author={Tardos, G\'{a}bor},
			title={Local chromatic number, Ky Fan's theorem and circular colorings},
			journal={Combinatorica},
			volume={26},
			date={2006},
			number={5},
			pages={587--626},
			issn={0209-9683},
			review={\MR{2279672}},
			doi={10.1007/s00493-006-0034-x},
		}

		
	\end{thebibliography}
\end{document}